\documentclass[10pt]{amsart}

\usepackage{amsmath,amsthm,amsfonts,amssymb}
\usepackage{graphicx}
\usepackage{hyperref}

\newtheorem{theorem}{Theorem}[section]
\newtheorem{lemma}[theorem]{Lemma}
\theoremstyle{definition}
\newtheorem{algorithm}[theorem]{Algorithm}
\newtheorem{corollary}[theorem]{Corollary}

\newtheorem{definition}[theorem]{Definition}

\newtheorem{remark}[theorem]{Remark}
\newtheorem{proposition}[theorem]{Proposition}

\newcommand{\rb}[1]{{\left( #1 \right)}}

\newcommand{\rbb}[1]{{\Bigl( #1 \Bigr)}}

\newcommand{\modb}[1]{{\left| #1 \right|}}

\newcommand{\modbb}[1]{{\Bigl| #1 \Bigr|}}

\newcommand{\cb}[1]{{\left\{ #1 \right\}}}

\newcommand{\cbb}[1]{{\Bigl\{ #1 \Bigr\}}}

\newcommand{\gp}[1]{{\left\langle #1 \right\rangle}}
\newcommand{\sqb}[1]{{\left\{ #1 \right\}}}

\def\MS{{\mathbb{S}}}

\def\CF{{\mathcal F}}

\def\ovx{{\overline{x}}}
\def\ovzeta{{\overline{\zeta}}}

\def\ovX{{\overline{X}}}
\def\ovv{{\overline{v}}}
\def\ovu{{\overline{u}}}

\def\ovR{{\overline{R}}}
\def\ovo{{\overline{0}}}

\def\ME{{\mathbb{E}}}
\def\MN{{\mathbb{N}}}
\def\MZ{{\mathbb{Z}}}
\def\MR{{\mathbb{R}}}
\def\BP{{\bf{P}}}

\def\IND{{\mathbf{1}}}

\title[SSLN for Graphs and Groups]{Strong law of large numbers on graphs and groups}

\date{\today}

\begin{document}

\author[N. Mosina]{Natalia Mosina}
\address{Department of Mathematics, CUNY/LAGCC, Long Island City, NY, USA}
\email{nmosina@lagcc.cuny.edu; mosina@math.columbia.edu}
\thanks{
The work of the first author started in Columbia University and was partially supported by the NSF grant DMS-06-01774.
The work continued in CUNY under support of the PSC-CUNY Grant Award 60014-40 41.
The work of the second author was partially supported by the NSF grant DMS-0914773.}

\author[A. Ushakov]{Alexander Ushakov}
\address{Department of Mathematics, Stevens Institute of Technology, Hoboken, NJ, USA}
\email{sasha.ushakov@gmail.com}

\maketitle

\begin{abstract}
We consider (graph-)group-valued random element $\xi$, discuss the
properties of a mean-set $\ME(\xi)$, and prove the
generalization of the strong law of large numbers for graphs and
groups.
Furthermore, we prove an analogue of the classical Chebyshev's
inequality for $\xi$ and Chernoff-like asymptotic bounds.
In addition,
we prove several results about configurations of mean-sets in
graphs and discuss computational problems together with methods
of computing mean-sets in practice and propose an algorithm for
such computation.

\end{abstract}




\section{Introduction}

Random objects with values in groups and graphs are constantly dealt
with in many areas of mathematics and theoretical computer science.
In particular, such objects are very important in group-based
cryptography (see \cite{MSU_book} or \cite{Dehornoy_survey} for
introduction to the subject). Having the notion of the average for
random group elements, generalized laws of large numbers for groups
with respect to this average together with results on the rate of
convergence in these laws would broaden the range of applications of
random group objects from both theoretical and practical point of
view. With a continuing development of group-based cryptography,
availability of such tools for analysis of probability measures and
their characteristics on groups becomes especially important. In
this paper, we develop these probabilistic tools for
finitely generated groups and propose practical algorithms for
computing mean values (or expectations) of group/graph-valued random
elements. The results of this paper form a new mathematical
framework for group-based cryptography and find applications to
security analysis of Sibert type authentication protocols (\cite{MosUsh:MeanAttack}).

The classical strong law of large numbers (SLLN) states that for independent
and identically distributed (i.i.d.) real-valued random variables
$\{\xi_i\}_{i=1}^\infty$
\begin{equation}\label{eq:slln}
\frac{1}{n}\sum_{i=1}^{n}\xi_i\rightarrow \ME(\xi_1)
\end{equation}
almost surely (with probability one) as $n\rightarrow\infty$,
provided that expectation $\ME(\xi_1)$ is finite (see
\cite{Billingsley}, \cite{Feller1971:IntroProbabilityTheory}, or
\cite{Skorohod:2004}). It is natural to pose a question about the
existence of counterparts of this result for different topological
spaces and/or algebraic structures, including groups. Starting from
the middle of the last century, there has been ongoing research,
following different trends, concerning the existence of such
generalizations of the SLLN. One line of this research investigates
random walks on groups (see Section \ref{se: slln_rw_on_groups} for
a brief list of relevant literature sources). The present work
follows another direction of that research -- the one which is
concerned with the problem of averaging in arbitrary metric spaces.
We generalize classical probability results to groups starting with
the concept of expectation (mean value) for group elements. Then we
prove the almost sure (with probability one) convergence, in some
appropriate sense, of sample (empirical) means for group/graph
random elements to the actual (theoretical) mean, thus, generalizing
the classical law (\ref{eq:slln}) and preserving its fundamental
idea. We supplement our results with the analogues of
Chebyshev and Chernoff-like bounds on the rate of convergence in the
SLLN for random graph/group elements.

\subsection{Historical Background}

Below we give a brief account of some developments concerning
probabilities and mean-values for various spaces as well as some
already existing generalizations of the strong law of large numbers
in order to highlight several stages of research that preceded our
work.
The reader willing to proceed to the core of our work right away may
skip this section and move on to Section \ref{se:Intro_corework}.

\subsubsection{Linear spaces} In 1935, Kolmogorov
\cite{Kolmogorov:1935} proposed to study probabilities in Banach
spaces. Later, the interpretation of stochastic processes as random
elements in certain function spaces inspired the study of laws of
large numbers for random variables taking values in linear
topological spaces. Banach spaces fit naturally into the context of
the strong law of large numbers because the average of $n$ elements
$x_1,\ldots,x_n$ in a Banach space is defined as
$n^{-1}(x_1+\ldots+x_n)$. In addition, Banach space provides
convergence, and Gelfand--Pettis integration provides the notion of
expectation of a random element (see \cite{Gelfand:1936} and
\cite{Pettis:1938}). It goes as follows. Let $X$ be a linear space
with norm $\|\cdot\|:X\rightarrow \MR$ and $X^\ast$ is the
topological dual of $X$. A random $X$-element $\xi$ is said to have
the expected value $\ME(\xi) \in X$  if
    $$\ME(f(\xi)) = f(\ME(\xi))$$
for every $f\in X^\ast$. Let $\{\xi_i\}_{i=1}^\infty$ be a sequence
of random $X$-elements. Without loss of generality, we may assume
that $\ME \xi_i = 0$ for every $i$. The strong law of large numbers
in a separable Banach space $X$ for a sequence of i.i.d. random
$X$-elements $\{\xi_i\}_{i=1}^\infty$ is first proved in
\cite{Mourier:1953}. It states that
    $$\lim_{n\rightarrow \infty} \|n^{-1}(\xi_1+\ldots+\xi_n)\| = \ME(\xi_1) = 0$$
with probability one.
The strong law of large numbers for i.i.d.
random elements in a Fr\'echet space was proved in
\cite{Ahmad:1965}.
A few other works discussing generalizations of
the strong law of large numbers in linear spaces are
\cite{Beck:1963}, \cite{Beck_Giesy:1970},
\cite{Taylor:1972}.

\subsubsection{Metric spaces} \label{se:FrechetMeanAndMore}

Unlike in linear spaces, in a general (non-linear) topological space
$X$, one has to do find some other ways to introduce the concept of
averaging and expectation. In $1948$, Fr\'{e}chet, \cite{Frechet48},
proposed to study probability theory in general metric spaces and
introduced a notion of a mean (sometimes called \textit{Fr\'{e}chet
mean}) of a probability measure \textit{$\mu$} on a complete metric
space $(X, d)$ as the minimizer of $E d^2(x,y)$, where
$$E d^2(x,y) = \int_X d^2(x,y)\mu(dy)$$
when it exists and is unique. If the minimizer is not unique, then
the set of minimizers can be considered.
If $\xi:\Omega \rightarrow X$ on a given probability space
$(\Omega,\CF,\BP)$ (see \cite{Billingsley},
\cite{Feller1971:IntroProbabilityTheory}) is a random element in
$(X, d)$ and if for some $\textbf{x} \in X$,
\begin{equation} \label{eq:IntroFrechetMeanGeneral}
E d^2(\xi,\textbf{x}) = \inf_{y \in X}E d^2(\xi,y) < \infty,
\end{equation}
then $\textbf{x}$ is called an {\em expected element} of $X$. These
generalizations were not met with much enthusiasm at the time (see
historical remarks on probabilities in infinite dimensional vector
spaces in \cite{Grenander}), and Fr\'{e}chet's suggestions, due to
the luck of their applications, underwent rather slow developments
in the middle of the last century.

Let us briefly mention some existing works on generalizing the
classical SLLN to a metric space $X$.
Let $\{\xi_i\}_{i=1}^\infty$ be a sequence of i.i.d. random elements
with values in $X$. Let the
expectation be defined as in
(\ref{eq:IntroFrechetMeanGeneral}), written as a set
$$\ME(X) = \displaystyle \cb{x\in X ~\bigg{|}~ E d^2(\xi,x) =
\inf_{y\in X}E d^2(\xi,y)}.$$
Define an empirical mean (average) of
elements $\xi_1(\omega),\ldots,\xi_n(\omega)$ to be the set
    $${\bf M}(\xi_1,\ldots,\xi_n) = \cb{x\in X ~\bigg{|}~  \sum_{i=1}^n d^2(x,\xi_i(\omega)) = \inf_{y\in X} \sum_{i=1}^n d^2(y,\xi_i(\omega)) }$$
One of the first works on
generalization of the SLLN for metric spaces is given in $1977$ by
Ziezold (\cite{Ziezold:1977}).
Ziezold considers a separable quasi-metric space $X$ with a finite
quasi-metric $d$. For a sequence of i.i.d. random $X$-elements
$\{\xi_i\}_{i=1}^\infty$ such that $E d^2(\xi,x)$ is finite for at
least one $x\in X$, he proves that inclusion
\begin{equation} \label{eq:ZiezoldInclusion}
 {\bf M (\omega)} = \bigcap_{k=1}^\infty \overline{\bigcup_{n=k}^\infty {\bf M}(\xi_1(\omega),\ldots,\xi_n(\omega))} \subseteq {\bf E(\xi_1)}
 \end{equation}
holds with probability one. Here, $\overline{\bigcup_{n=k}^\infty
{\bf M}(\xi_1(\omega),\ldots,\xi_n(\omega))}$ is the closure of the
union of ${\bf M}(\xi_1(\omega),\ldots,\xi_n(\omega))$'s. He also
shows that, in general, the equality does not hold (for a finite
quasi-metric space). In 1981, Sverdrup-Thygeson
(\cite{Sverdrup-Thygeson:1981}) proves inclusion
(\ref{eq:ZiezoldInclusion}) for compact connected metric spaces and
shows that the equality does not hold in general (for a metric
space) when the minimizer in (\ref{eq:IntroFrechetMeanGeneral}) is
not unique.
In 2003, Bhattacharya and Patrangenaru in \cite[Theorem
2.3]{Bhattacharya_Patrangenaru:2003} prove equality in
(\ref{eq:ZiezoldInclusion}) for the unique minimizer in
(\ref{eq:IntroFrechetMeanGeneral}) for metric spaces $X$ such that
every closed bounded subset of $X$ is compact, improving Ziezold's
and Sverdrup-Thygeson's results.

\vspace{2mm}
\noindent{\bf Manifolds.} As the need for statistical analysis for
spaces with differential geometric structure was arising,
statistical inference on Riemannian manifolds started to develop
rapidly, especially due to applications in statistical theory of
shapes and image analysis. These applications evolve around the
concept of averaging. See \cite{Kendall:book1999} for an
introduction into shape theory. The interested reader may also refer
to
\cite{Huckemann:2010},
for instance.

There are two main approaches to averaging of elements on a
manifold. Every Riemannian manifold $X$ is a metric space and hence
one can use constructions from the previous section to define the
notion of a mean. Fr\'{e}chet mean of a probability measure on a
manifold is also known as an {\em intrinsic mean}
\cite{Bhattacharya_Patrangenaru:2003}. Non-uniqueness of the
intrinsic mean is a source of different technical problems. Also,
the intrinsic mean, even when unique, is often very difficult to
compute in practice.

On the other hand, a manifold $X$ can also be looked at as a
submanifold of some Euclidean space $\MR^k$ and one can define a
mean relative to this inclusion. Let $\tau:X\rightarrow \MR^k$ be an
embedding of $X$ into Euclidean space $(\MR^k,d_0)$. A point $p \in
\MR^k$ is called {\em nonfocal} if there exists a unique $x\in
\tau(X)$ such that
    $d_0(p,x) = d_0(p,\tau(X)).$
Let $\mu$ be a probability measure on $X$, $\mu'$ a probability
measure on $\MR^k$ induced by $\tau$, and $x^\ast \in \MR^k$ the
expectation of $\mu'$. We say that the measure $\mu$ is nonfocal if
$x^\ast$ is a nonfocal point. For a nonfocal probability measure
$\mu$ on $X$ we define the mean as $\tau^{-1}(x^\ast)$. In
\cite{Bhattacharya_Patrangenaru:2003} the authors prove the strong
law of large numbers for the intrinsic and extrinsic means on
manifolds.

\subsubsection{K-means} \label{se:history_other_means}

A notion of a mean (or a mean-set) can be generalized into $k$-mean.
Let $B$ be a Banach space with a norm $\|\cdot\|$ and $k\in\MN$. For
a set $H = \{h_1,\ldots,h_k\}$ we define a partition of $B$ as
follows
    $$S_i = \{x\in B \mid \|x-h_i\|\le \|x-h_j\| \mbox{ for every } j=1,\ldots,k\}\setminus(S_1\cup\ldots\cup S_{i-1})$$
where $i=1,\ldots,k$ and a function $\pi_{H}:B\rightarrow B$
    $$\pi_{H}(x) = \sum_{i=1}^k h_i \cdot \IND_{S_i}(x)$$
where $\IND_{S_i}$ is the indicator function of $S_i$. Fix a suitable
non-decreasing function $\Phi:\MR^+\rightarrow \MR^+$ (e.g.,
$\Phi(x)=x$ or $\Phi(x)=x^2$) and for a probability measure $\mu$ on
$B$ define a number
    $$M(H) = \int_{x\in B} \Phi(\|x-\pi_H(x)\|) d\mu(x).$$
A set $H_0$ of $k$ elements that minimizes the value of $M$ is
called a {\em $k$-mean} of a probability distribution $\mu$. In
general, there can be several minimizers, which leads to technical
complications. In 1988, Cuesta and Matran (\cite{
Cuesta_Matran:1988}) proved that empirical $k$-means converge to the
$k$-mean $H_0$ of $\mu$ under the assumption that the $k$-mean is
unique.

It is straightforward to generalize a notion of a $k$-mean to a
general metric space $(X,d)$. Indeed, if we put
    $$S_i = \{x\in B \mid d(x,h_i)\le d(x,h_j) \mbox{ for every } j=1,\ldots,k\}\setminus(S_1\cup\ldots\cup S_{i-1})$$
and
    $$M(H) = \int_{x\in B} \Phi(d(x,\pi_H(x))) d\mu(x),$$
then we get a similar notion. This type of $k$-means, with
$\Phi(x)=x$, was considered by Rubinshtein in 1995 in
\cite{Rubinshtein1995} where it was called the $k$-center.

As we can see, in general, depending on the research goals, one can
define mean values on a given metric space $(X,d)$ using any powers
of $d$, i.e., instead of dealing with minimization of $E d^2(\xi,x)$
in (\ref{eq:IntroFrechetMeanGeneral}), one can work with a very
similar functional by minimizing $E d^r(\xi,x)$ for any $r>0$ if
necessary.

\subsubsection{Probabilities on algebraic structures}

Metrics and probabilities on algebraic structures have been studied
from different perspectives. One source to look at is the book of M.
Gromov \cite{Gromov_book:1999}. The reader can find some applications
of Fr\'{e}chet mean in statistical analysis of partially ranked data
(such as elements of symmetric groups and homogeneous spaces) in the
book of Diaconis (\cite{Diaconis88}).
An extensive historical background of the studies of probabilities
on algebraic structures is given in \cite{Grenander}, where the
author considers probabilities for stochastic semi-groups, compact
and commutative stochastic groups, stochastic Lie groups, and
locally compact stochastic groups employing the techniques of
Fourier analysis to obtain limit theorems for convolutions of
probability distributions. The reader interested in the question of
defining probabilities on groups can find several approaches to this
issue in \cite{BoMS}.

\subsubsection{Random walks on groups} \label{se: slln_rw_on_groups}

One way to generalize the strong law of large numbers for groups is
to study the asymptotic behavior of the products $g_1 g_2 \ldots
g_n$, where $\{g_i\}_{i=1}^\infty$ is a sequence of i.i.d. random
group elements, the so-called random walk on a group. The reader can
consult \cite{Woess:1994} for an introduction to random walks on
groups. In 1960, Furstenberg and Kesten
(\cite{Furstenberg_Kesten:1960}) prove the generalization of the
strong law of large numbers for random matrices. Namely, they show
that the limit
    $$\lim_{n\rightarrow\infty}\frac{\log ||g_1 g_2 \ldots g_n||}{n}$$
exists with probability one, with some restrictive conditions on the
entries of $g_i$, without computing the limit explicitly. In 1963,
Furstenberg solved this problem for normalized products of random
matrices in terms of stationary measures (\cite{Furstenberg:1963}).
Computational techniques that would allow to compute these measures
are investigated in \cite{Pincus:1985}. A concise account of a
number of illuminating results in the direction of the
generalization of the SLLN to groups can be found in
\cite{Karlsson-Ledrappier:2006}, where the authors prove the theorem
about the directional distribution of the product $g_1 g_2 \ldots
g_n$, thus, proving a general law of large numbers for random walks
on general groups. The authors call it a multiplicative ergodic
theorem or a general, noncommutative law of large numbers (see
\cite{Karlsson-Ledrappier:2006} for the precise statement).

\subsection{The core of our work}\label{se:Intro_corework}

Motivated by applications to group-based cryptanalysis, we study
Fr\'{e}chet type mean values and their properties in graph/group
theoretic settings.

Let $\Gamma = (V(\Gamma), E(\Gamma))$ be a locally finite graph and
$(\Omega,\CF,\BP)$ a given probability space. A random
$\Gamma$-element $\xi$ is a measurable function
$\xi:\Omega\rightarrow V(\Gamma)$. This random $\Gamma$-element
$\xi$ induces an atomic probability measure
$\mu:V(\Gamma)\rightarrow [0,1]$ on $V(\Gamma)$ in a usual way:
    $$\mu(v) = \BP(\{\omega\in \Omega \mid \xi(\omega)=v \}), ~ v\in V(\Gamma).$$
Next, we introduce a {\em weight function} $M_\xi:V(\Gamma)
\rightarrow \MR$ by
   $$M_\xi(v) = \ME d^2(v,\xi) = \sum_{s\in V(\Gamma)} d^2(v,s) \mu(s),$$
where $d(v,s)$ is the distance between $v$ and $s$ in $\Gamma$, and
note that, trivially, the domain of definition of $M_\xi(\cdot)$ is
either the whole $V(\Gamma)$ (in which case we say that $M$ is {\em
totally defined}) or $\emptyset$. The domain of $M$ is the set
    $$domain(M) = \cb{v\in\ V(\Gamma) ~\bigg{|}~ \sum_{s\in V(\Gamma)} d^2(v,s) \mu_\xi(s)<\infty }.$$
In the case when $domain(M_\xi) = V(\Gamma)$, we define the {\em
mean-set} of $\xi$ to be
\begin{equation}\label{eq:mean-set}
    \ME (\xi) = \{v \in V(\Gamma) \mid  M_\xi(v) \le M_\xi(u), ~~ \forall u\in V(\Gamma) \}.
\end{equation}
The above definition of $\ME (\xi)$,
$\xi:\Omega\rightarrow V(\Gamma)$, provides the corresponding notion
of a mean (average, expectation) for finitely generated groups via
their Cayley graphs.

Once we have the notion of mean-set for group-valued random
elements, we notice that it satisfies the so-called ``shift''
property; namely,
\begin{equation}\label{eq:shift_prop}
    \ME (g\xi) = g\ME (\xi), \forall g \in G
\end{equation}
which is analogous to the linearity property of a classical
expectation for real-valued random variables.

Next, for a sample $\xi_1(\omega), \ldots, \xi_n(\omega)$ of i.i.d.
random $\Gamma$-elements we define a {\em relative frequency}
$\mu_n(u;\omega)= \mu_n(u)$ with which the value $u\in V(\Gamma)$
occurs in the sample above:
    $$\mu_n(u) = \frac{1}{n}|\{i \mid \xi_i=u\}|.$$
Relative frequency $\mu_n$ is a probability measure on $\Gamma$, and
we can define {\em empirical (sampling) weight function} (random
weight) as
   $$M_n(v) = \sum_{s\in V(\Gamma)} d^2(v,s) \mu_n(s).$$
Going further, we define an {\em empirical mean} or {\em sample
mean-set} of the sample $\xi_1,\ldots,\xi_n$ to be the set of
vertices
    $$\MS(\xi_1,\ldots,\xi_n) = \{v \in V(\Gamma) \mid  M_n(v) \le M_n(u), ~~ \forall u\in V(\Gamma)\}.$$
The function $\MS(\xi_1,\ldots,\xi_n)$ on graphs is an analogue of
the average function $(x_1,\ldots,x_n)\mapsto(x_1+\ldots+x_n)/n$ for
$x_1,\ldots,x_n \in \MR$. We let $\MS_n = \MS(\xi_1,\ldots,\xi_n)$.
With these notions at hand, we first formulate and prove the
following generalization of the strong law of large numbers for
graphs and groups with one-point mean sets.

\vspace{2mm}\noindent{\bf Theorem A.} {\bf(Strong Law of Large
Numbers for graphs )} {\em Let $\Gamma$ be a locally-finite connected graph and
$\{\xi_i\}_{i=1}^\infty$ a sequence of i.i.d. random
$\Gamma$-elements. If the weight function $M_{\xi_1}(\cdot)$ is totally defined and
    $\ME(\xi_1) = \{v\}$
for some $v\in V(\Gamma)$, then
    $$\lim_{n\rightarrow\infty} \MS(\xi_1,\ldots,\xi_n) = \ME(\xi_1)$$
with probability one.
}

Next, we improve this result and prove the generalized law of large
numbers for groups for the case when $|\ME\xi|>1$ (see Section
\ref{slln2or3vert}).
The simplest version of multi-vertex SLLN in terms of {\em limsup}
is as follows:

\vspace{2mm}\noindent{\bf Theorem B.} {\bf(Multi-vertex SLLN for
graphs )}
{\em Let $\Gamma$ be a locally-finite connected graph and
$\{\xi_i\}_{i=1}^\infty$ a sequence of i.i.d.
random $\Gamma$-elements. Assume that the
weight function $M$ is totally defined and $\ME(\xi) =
\{v_1,\ldots,v_k\}$, where $k\ge 4$. If the random walk $\ovR^1$
associated to $v_1$ is genuinely $(k-1)$-dimensional, then
    $$\limsup_{n\rightarrow \infty}\MS_n = \ME(\xi_1)$$
holds with probability one.}

In addition, we prove analogues of classical Chebyshev's inequality
and Chernoff-like bounds for a graph-(group-) random element $\xi$.

\vspace{2mm}\noindent{\bf Theorem C.}  {\bf(Chebyshev's inequality
for graphs)} {\em Let $\Gamma$ be a locally-finite connected graph
and $\{\xi_i\}_{i=1}^\infty$ a sequence of i.i.d.
random $\Gamma$-elements. If the weight
function $M_{\xi_1}(\cdot)$ is totally defined then there exists a
constant $C = C(\Gamma,\xi_1)>0$ such that
    $$\displaystyle \BP\rbb{ \MS(\xi_1,\ldots,\xi_n) \not\subseteq \ME(\xi) } \le \frac{C}{n}.$$
}

\vspace{2mm}\noindent{\bf Theorem D.}  {\bf(Chernoff-like bounds for
graphs)} {\em Let $\Gamma$ be a locally-finite connected graph and
$\{\xi_i\}_{i=1}^\infty$ a sequence of i.i.d. random
$\Gamma$-elements. If the weight function
$M_{\xi_1}(\cdot)$ is totally defined and $\mu_{\xi_1}$ has finite
support, then for some constant $C>0$
    $$\BP\rbb{ \MS(\xi_1,\ldots, \xi_n) \not\subseteq \ME(\xi) } \le O(e^{-Cn}).$$
}

\subsection{Outline}
In Section \ref{se:mean}, we give basic definitions and discuss some
properties of the newly defined objects. Next, we turn to the
formulation and the proof of the strong law of large numbers on
graphs and groups. These tasks are carried out in Section
\ref{slln}. Chebyshev's inequality and Chernoff-like bounds for
graphs are proved in Section \ref{se:Chebyshev}. In Section
\ref{se:conf_centers}, we consider configurations of mean-sets in
graphs and their applications to trees and free groups.
Section \ref{se: Computation_meanset} deals with computational
problems and methods of computing $\ME(\xi)$. In particular, we
propose an algorithm and prove that this algorithm finds a central
point for trees.
Finally, in Section \ref{se: experiments}
we perform series
of experiments in which we compute the sample mean-sets of randomly
generated samples of $n$ random elements and observe the convergence
of the sample mean-set to the actual mean.

\section{Mean (expectation) of a group-valued random element}
\label{se:mean}

Let $(\Omega,\CF,\BP)$ be a given probability space and $G$ a finitely
generated group.
In this section, we define the notion of expectation for
graph random elements in the sense of Fr\'{e}chet mean set, which is
one of the possible ways to look at mean values (see
Section \ref{se:FrechetMeanAndMore}). The same definition will hold
for random elements $\xi:\Omega \rightarrow G$ on groups.
We also discuss properties of the mean sets on groups.
In particular, we prove
that for our expectation $\ME$, we have
    $\ME(g\xi) = g\ME(\xi).$


\subsection{The mean set in a graph}

\label{se:Graph_mean}

Let $\Gamma = (V(\Gamma), E(\Gamma))$ be a locally finite connected graph.
A random
$\Gamma$-element $\xi$ is a measurable function
$\xi:\Omega\rightarrow V(\Gamma)$. The random element
$\xi$ induces an atomic probability measure
$\mu:V(\Gamma)\rightarrow [0,1]$ on $V(\Gamma)$ in a usual way:
\begin{equation}\label{eq:measure_mu}
\mu(v) = \mu_\xi(v) = \BP(\{\omega\in \Omega \mid \xi(\omega)=v \}), ~ v\in V(\Gamma).
\end{equation}
Next, we introduce a {\em weight function} $M_\xi:V(\Gamma)
\rightarrow \MR$ by
   $$M_\xi(v) = \ME d^2(v,\xi) = \sum_{s\in V(\Gamma)} d^2(v,s) \mu(s),$$
where $d(v,s)$ is the distance between $v$ and $s$ in $\Gamma$.
If $M_{\xi}(v)$ is finite, then we say that the {\em weight function}
$M_{\xi}$ is defined at $v$.
The domain of $M$ is the set
    $$domain(M) = \cb{v\in\ V(\Gamma) ~\bigg{|}~ \sum_{s\in V(\Gamma)} d^2(v,s) \mu_\xi(s)<\infty }.$$
The case of interest of course, is when $M_{\xi}(v)$ is {\em totally
defined}, meaning that $M_{\xi}(v)$ is finite at every $v\in V(\Gamma)$.

\begin{definition}\label{def:meanset}
Let $\xi$ be a random $\Gamma$-element such that $M_{\xi}(\cdot)$ is totally defined.
The set of vertices $v\in \Gamma$ that minimize the value of $M_\xi$
\begin{equation}\label{eq:meanset}
   \ME (\xi) = \{v \in V(\Gamma) \mid M_{\xi}(v) \le M_{\xi}(u),  ~~  \forall u\in V(\Gamma)\},
\end{equation}
is called the {\em mean-set} (or the {\em center-set}, or {\em average}) of $\xi$.
\end{definition}

Very often we leave the random element $\xi$ in the background to
shorten the notation and write $M(v)$ instead of $M_{\xi}(v)$.
Moreover, we write $\ME (\mu)$ instead of $\ME (\xi)$ sometimes and
speak of the mean set of distribution $\mu$ induced by $\xi$ on
$V(\Gamma)$.


\begin{lemma}\label{le:p_defined}
Let $\Gamma$ be a connected graph, $\xi$ a random $\Gamma$-element,
and $u,v$ adjacent vertices in $\Gamma$.
If $M(u)<\infty$, then $M(v)<\infty$.
\end{lemma}

\begin{proof}
Easily follows from the definition of $M$ and the triangle inequality.
\end{proof}

\begin{corollary}\label{co:M_domain}
Let $\Gamma$ be a connected graph and $\xi$ a random $\Gamma$-element.
Then either $domain(M) = V(\Gamma)$ or $domain(M) = \emptyset$.
\end{corollary}

\begin{lemma}\label{le:E_finite}
Let $\Gamma$ be a connected locally finite graph and
$\xi$ a random $\Gamma$-element. If $M_\xi$ is totally defined,
then $0<|\ME(\xi)|<\infty$.
\end{lemma}

\begin{proof}
Let $\mu$ be a measure of (\ref{eq:measure_mu}) induced on $\Gamma$
by $\xi$. For an arbitrary but fixed vertex $v \in \Gamma$, the
weight function
    $$M(v) = \sum_{i\in V(\Gamma)} d^2(v,i) \mu(i) = \sum_{n=0}^\infty \rb{ n^2 \sum_{i\in V(\Gamma), d(v,i)=n} \mu(i)}$$
is defined at $v$ by assumption. Choose $r \in \MN$ such that
    $$\frac{1}{2} M(v) \le \sum_{n=0}^r \rb{ n^2 \sum_{i\in V(\Gamma), d(v,i)=n} \mu(i) } = \sum_{i\in B_v(r)} d^2(v,i) \mu(i),$$
where
\begin{equation}  \label{eq:ball}
   B_v(r) = \{i\in V(\Gamma) \mid d(v,i)\le r \}
\end{equation}
is the {\em ball} in $\Gamma$ of radius $r$ centered at $v$. If we
take a vertex $u$ such that $d(u,v) \ge 3r$, then using the triangle
inequality, we obtain the following lower bound:
    $$M(u) = \sum_{i\in V(\Gamma)} d^2(u,i) \mu(i)\ge \sum_{i\in B_v(r)} [2r]^2 \mu(i) + \sum_{i\not\in B_v(r)} d^2(u,i) \mu(i)
    \ge 4\sum_{i\in B_v(r)} d^2(v,i) \mu(i) \ge 2 M(v).$$
Thus, $d(v,u) \ge 3r$ implies $u\not \in \ME (\xi)$ and, hence,
$\ME(\xi) \subseteq B_v(3r)$. Since the graph $\Gamma$ is locally
finite, it follows that the sets $B_v(3r)$ and $\ME (\xi)$ are
finite. This implies that the function $M$ attains its minimal value
in $B_v(3r)$ and hence $\ME(\xi) \ne \emptyset$.
\end{proof}

\subsection{The mean set in a group}
\label{se:Group_mean}

Let $G$ be a group and $X \subseteq G$ a finite generating set for $G$.
The choice of $X$ naturally determines a distance $d_X$ on $G$
via its Cayley graph $C_G(X)$.
Hence Definition  \ref{def:meanset} gives us a notion of a mean set for
a random $G$-element.
It follows from the definition of the distance $d_X$ that for
any $a,b,g \in G$ the equality
\begin{equation}\label{eq:distance}
    d_X(a,b) = d_X(ga,gb)
\end{equation}
holds. This equality implies that $\ME(\xi)$ possesses the desirable
property $\ME(g\xi)=g\ME(\xi)$, as the following proposition shows.

\begin{proposition}[``Shift'' Property]      \label{pr:g_shift}
{\em Let $G$ be a group and $g \in G$. Suppose that
$(\Omega,\CF,\BP)$ is a given probability space and $\xi:\Omega
\rightarrow G$ a $G$-valued random element on $\Omega$. Then for the
random element $\xi_g$ defined by $\xi_g(\omega) = g\xi(\omega)$ we
have
    $\ME (\xi_g) = g\ME (\xi).$
}
\end{proposition}

\begin{proof}
Let $\mu_{\xi_g}$ be the measure induced on $G$ by $\xi_g$, in the
manner of (\ref{eq:measure_mu}). It follows from the definition of
$\xi_g$ that for any $h\in G$
    $$\mu_{\xi_g}(h) = \BP(\{\omega \mid \xi_g(\omega)=h \}) = \BP(\{\omega \mid g\xi(\omega)=h \})
    = \BP(\{\omega \mid \xi(\omega)=g^{-1}h \}) = \mu_\xi(g^{-1}h).$$
This, together with (\ref{eq:distance}), implies that for any $h \in
G$
    $$M_{\xi_g}(h) = \sum_{i\in G} d^2(h,i) \mu_{\xi_g}(i) = \sum_{i\in G} d^2(g^{-1}h,g^{-1}i) \mu_{\xi}(g^{-1}i) = \sum_{i\in G} d^2(g^{-1}h,i) \mu_{\xi}(i) = M_{\xi}(g^{-1}h).$$
Hence, the equality $M_{\xi_g}(h) = M_{\xi}(g^{-1}h)$ holds for any
random element $\xi$ and $g,h \in G$. Therefore, for any
$h_1, h_2 \in G$,
    $M_{\xi_g}(h_1) < M_{\xi_g}(h_2) \Leftrightarrow M_{\xi}(g^{-1}h_1) < M_{\xi}(g^{-1}h_2)$
and
    $$\ME (\xi_g) = \cbb{h\in G \mid  M_{\xi_g}(h)\le M_{\xi_g}(f),~~ \forall f\in G } = \cbb{h\in G \mid  M_{\xi}(g^{-1}h)\le M_{\xi}(g^{-1}f), ~~\forall f\in G } = $$
    $$= \cbb{h\in G \mid M_{\xi}(g^{-1}h)\le M_{\xi}(f),~~ \forall f\in G } = \cbb{gh\in G \mid M_{\xi}(h)\le M_{\xi}(f),~~ \forall f\in G } = g\ME (\xi).$$
\end{proof}

The equality $d_X(a,b) = d_X(ag,bg)$ does not hold for a general
group $G = \gp{X}$. It holds for abelian groups.

\begin{proposition}\label{pr:g_shift2}
{\em Let $G$ be an abelian group and $g \in G$. Suppose that
$(\Omega,\CF,\BP)$ is a probability space and $\xi:\Omega
\rightarrow G$ a $G$-valued random element on $\Omega$. Then for the
random element $\xi_g$ defined by $\xi_g(\omega) = \xi(\omega)g$ we
have
    $\ME (\xi_g) = (\ME (\xi)) g.$
}
\end{proposition}

\subsection{Other possible definitions of $\ME$} \label{se: other_means}

There are other possible definitions of $\ME$ for which the
statement of Proposition \ref{pr:g_shift} (and other results of Section \ref{slln}) holds. Let $c$ be a
positive integer. By analogy to the function $M_\xi(v)$, define a
{\em weight function} $M_\xi^{(c)}(v)$ of class $c$ by
    $$M_\xi^{(c)}(v) = \sum_{i\in V(\Gamma)} d^c(v,i) \mu(i)$$
and the mean-set $\ME^{(c)} (\xi)$ of class $c$ to be
    $$\ME^{(c)} (\xi) = \{v \in V(\Gamma) \mid  M^{(c)}(v) \le M^{(c)}(u), ~~ \forall u\in V(\Gamma) \}.$$
It is straightforward to check that all the statements of the
previous section hold for $M_\xi^{(c)}(\cdot)$ and $\ME^{(c)}
(\xi)$. In fact, it is not hard to see that when $c=1$, we have a
counterpart of the median of the distribution $\mu$.
Next proposition shows that our $\ME$ agrees with the classical
definition of the expectation on $\MZ$ in the following sense.

\begin{proposition}\label{pr:connection_classical_mean}
{\em Let $\xi:\Omega\rightarrow \MZ$ be an integer-valued random
variable with classical expectation
    $\mathfrak{m} = \sum_{n\in \MZ} n\BP(\xi=n).$
Assume that $M\equiv M^{(2)}_\xi$ is defined on $\MZ$. Then
$1\le|\ME^{(2)} \xi| \le 2$ and for any $v\in \ME^{(2)} (\xi)$, we
have $|\mathfrak{m}-v| \le\frac{1}{2}$. }
\end{proposition}

\begin{proof}
Straightforward.
\end{proof}

\begin{remark}
Observe that $\ME^{(2)}$ does not coincide with the classical mean
in $\MR^2$. Recall that the classical mean in $\MR^2$ is defined
coordinate-wise, i.e., the mean of $(x_1,y_1),\ldots,(x_n,y_n)$ is a
point in $\MR^2$ defined by $\rb{\ME X,\ME Y}.$ For example,
consider the distribution on $\MZ^2$ such that $\mu(0,0) = \mu(0,3)
= \mu(3,0) = 1/3$ and for all other points $\mu=0$. Then the
classical mean
is the point $(1,1)$, and the mean-set $\ME^{(2)}$ is the point
$(0,0)$.
\end{remark}


\section{Strong Law of Large Numbers} \label{slln}

Let $\xi_1, \ldots, \xi_n$ be a sample of independent and
identically distributed (i.i.d.) graph-valued random elements
$\xi_i:\Omega\rightarrow V(\Gamma)$ defined on a given probability
space $(\Omega,\CF,\BP)$. For every $\omega\in\Omega$, let
$\mu_n(u;\omega)$ be the {\em relative frequency}
\begin{equation} \label{eq:mu_n}
    \mu_n(u;\omega) = \frac{|\{i\mid \xi_i(\omega)=u,~~ 1\le i\le n\}|}{n}
\end{equation}
with which the value $u\in V(\Gamma)$ occurs in the random sample
$\xi_1(\omega), \ldots, \xi_n(\omega)$. We shall suppress the
argument $\omega\in\Omega$ to ease notation, and let
    $$M_{n}(v) = \sum_{i\in V(\Gamma)} d^2(v,i) \mu_n(i)$$
be the {\em sampling weight}, corresponding to $v\in V(\Gamma)$,
and $M_{n}(\cdot)$ the resulting {\em sampling weight function}.

\begin{definition}
The set of vertices
$$\MS_n = \MS(\xi_1,\ldots,\xi_n) = \{v \in V(\Gamma) \mid  M_n(v) \le M_n(u), ~~ \forall u\in V(\Gamma)\}$$
is called the {\em sample mean-set} (or {\em sample center-set}, or average)
of the vertices $\xi_1,\ldots,\xi_n$.
\end{definition}

\begin{lemma}  \label{le:convergence_M_n}
Let $\Gamma$ be a locally-finite connected graph, $v\in V(\Gamma)$,
and $\{\xi_i\}_{i=1}^\infty$ a sequence of i.i.d. random $\Gamma$-elements such that
$M_{\xi_1}(v)$ is defined. Then
\begin{equation}\label{convergence_fixed_v}
    \BP\rbb{M_n(v) \rightarrow M(v) \mbox{ as } n\rightarrow\infty } = 1.
\end{equation}
\end{lemma}
\begin{proof}
For every $v\in V(\Gamma)$, $M(v)$ is the expectation of the random
variable $d^2(v,\xi_1)$.
The result follows by the strong law of large numbers applied to
$\{d^2(v,\xi_i)\}_{i=1}^\infty$.
\end{proof}

It is important to notice that in general the convergence in Lemma
\ref{le:convergence_M_n} is not uniform in a sense that, for some
distribution $\mu$ on a locally finite (infinite) graph $\Gamma$ and
some $\varepsilon>0$, it is possible that
    $$\BP\rbb{\exists N \mbox{ s.t. } \forall n>N ~\forall v\in V(\Gamma),~~ |M_n(v) - M(v)|<\varepsilon } < 1.$$
In other words, the convergence for every vertex, as in Lemma
\ref{le:convergence_M_n}, is insufficient to prove the strong law of
large numbers, stated in introduction. Next lemma is a key tool in
the proof of our strong law of large numbers.

\begin{lemma} [Separation Lemma]   \label{pr:separation_M}
Let $\Gamma$ be a locally-finite connected graph and
$\{\xi_i\}_{i=1}^\infty$ a sequence of i.i.d. random
$\Gamma$-elements. If the weight function $M_{\xi_1}(\cdot)$ is totally defined, then
    $$\BP\rbb{\exists N \mbox{ s.t. }\forall n>N,~ \max_{v\in\ME(\xi_1)} M_n(v) < \inf_{u\in V(\Gamma)\setminus\ME(\xi_1)} M_n(u) }  = 1.$$
\end{lemma}

\begin{proof}
Our goal is to prove that for some $\delta>0$
\begin{equation}\label{eq:main_eq_singleton}
    \BP \rbb{\exists N~ \forall n>N~ \forall v\in\ME(\xi_1),~ \forall u\in V(\Gamma)\setminus\ME(\xi_1),~~ M_n(u)-M_n(v)\ge\delta} = 1.
\end{equation}
We prove the formula above in two stages. In the first stage we show
that for some fixed $v_0 \in \ME(\xi_1)$ and for sufficiently large
number $m>0$ the following holds
\begin{equation}\label{eq:main_eq_singleton_a}
    \BP \rbb{\exists N \mbox{ s.t. }  \forall n>N~ \forall v\in\ME(\xi_1),~\forall u\in V(\Gamma)\setminus B_{v_0}(m),~~ M_n(u)-M_n(v)\ge\delta} = 1
\end{equation}
in the notation of (\ref{eq:ball}). In the second stage we prove
that
\begin{equation}\label{eq:main_eq_singleton_b}
    \BP \rbb{\exists N \mbox{ s.t. }  \forall n>N~ \forall v\in\ME(\xi_1),~\forall u\in B_{v_0}(m)\setminus\ME(\xi_1),~~ M_n(u)-M_n(v)\ge\delta} = 1
\end{equation}
Having the formulae above proved we immediately deduce that
(\ref{eq:main_eq_singleton}) holds using $\sigma$-additivity of
measure.

Let $v_0\in \ME (\xi_1)$ and $\mu$ be the probability measure on
$\Gamma$ induced by $\xi_1$, as in (\ref{eq:measure_mu}).
Since the weight function $M(\cdot)$ is defined at $v_0$, we can
choose $r \in \MR$ as in Lemma \ref{le:E_finite}, such that
$\frac{1}{2} M(v_0) \le \sum_{i\in B_{v_0}(r)} d^2(v_0,i) \mu(i)$.
Put $m=3r$. In Lemma \ref{le:E_finite} we proved that, if a vertex
$u$ is such that $d(u,v_0) \ge 3r$, then
\begin{equation}\label{M_u_ge_2M_v}
    M(u) = \sum_{i\in V(\Gamma)} d^2(u,i) \mu(i) \ge 4\sum_{i\in B_{v_0}(r)} d^2(u,i) \mu(i) \ge 2 M(v_0).
\end{equation}
It implies that $\ME (\xi_1) \subseteq B_{v_0}(3r)$.

Since $\Gamma$ is locally finite, the set $B_{v_0}(r)$ of
(\ref{eq:ball}) is finite. We also know from the SLLN for the
relative frequencies $\mu_n(u)$ that $\mu_n(u)
\stackrel{a.s.}{\rightarrow} \mu(u)$ as $n\rightarrow\infty$. These
facts imply that for any $\varepsilon > 0$, the event
\begin{equation}\label{eq:mu_bound}
   C_\varepsilon = \sqb{\exists N=N(\varepsilon),~~ \forall n>N,~~ \forall u \in B_{v_0}(r),~~ |\mu_n(u)-\mu(u)| < \varepsilon}
\end{equation}
has probability one. In particular, this is true for
    $\varepsilon=\varepsilon^\ast=\frac{1}{4} \min \{\mu(u) \mid u \in B_{v_0}(r),~ \mu(u)\ne 0\},$
and the event $C_{\varepsilon^\ast}$ is a subset of
\begin{equation}\label{M_n_ge_3/2M}
    \sqb{\exists N=N(\varepsilon^\ast),~~ \forall n>N,~~ \forall u\in V(\Gamma) \setminus B_{v_0}(3r),~~ M_n(u) \ge \frac{3}{2}M(v_0)
    }.
\end{equation}
Indeed, on the event $C_{\varepsilon^\ast}$, as in
(\ref{eq:mu_bound}), we have  $\mu_n(i)\geq \frac{3}{4}\mu(i)$, $i
\in B_{v_0}(r)$. Using this fact together with (\ref{M_u_ge_2M_v}),
we can write
    $$M_n(u) = \sum_{i\in V(\Gamma)} d^2(u,i) \mu_n(i) \ge 4\sum_{i\in B_{v_0}(r)} d^2(u,i) \mu_n(i) \ge 3\sum_{i\in B_{v_0}(r)} d^2(u,i) \mu(i) \ge \frac{3}{2} M(v_0).$$
Thus we have
\begin{equation}\label{eq:main_eq_singleton56}
    \BP\rb{\exists N \mbox{ s.t. }  \forall n>N,~ \forall u\in V(\Gamma) \setminus B_{v_0}(3r),~~ M_n(u) \ge \frac{3}{2}M(v_0) } = 1.
\end{equation}
By Lemma \ref{le:convergence_M_n}, for any $v\in V(\Gamma)$ and any
$\varepsilon>0$, we have
    $$\BP\rbb{\exists N=N(\varepsilon),~~ \forall n>N,~~ |M_n(v)-M(v)| < \varepsilon} = 1$$
and, since $B_{v_0}(3r)$ is a finite set, we have simultaneous
convergence for all vertices in $B_{v_0}(3r)$, i.e.,
\begin{equation}\label{eq:conv_mv}
    \BP\rbb{\exists N=N(\varepsilon),~~ \forall n>N,~~ \forall v\in B_{v_0}(3r),~~ |M_n(v)-M(v)| < \varepsilon} = 1.
\end{equation}
In particular, remembering that $\ME(\xi_1) \subseteq B_{v_0}(3r)$,
for $\varepsilon = M(v_0)/4$,
\begin{equation}\label{eq:main_eq_singleton11}
    \BP\rb{\exists N=N(\varepsilon),~~ \forall n>N, ~\forall v\in \ME(\xi_1),~~ \frac{3}{4}M(v) < M_n(v) < \frac{5}{4}M(v) } = 1.
\end{equation}
Finally, we notice that on the intersection of the events in
(\ref{eq:main_eq_singleton56}) and (\ref{eq:main_eq_singleton11}),
we have
    $$M_n(u)-M_n(v) \ge \frac{3}{2}M(v) - \frac{5}{4}M(v) = \frac{1}{4}M(v)=\frac{1}{4}M(v_0),$$
by the virtue of the fact that $M(v_0)=M(v)$ (as both $v_0, v \in
\ME(\xi_1)$), and formula (\ref{eq:main_eq_singleton_a}) holds for
any $\delta$ such that $\delta \le \frac{1}{4}M(v_0)$.

For the second part of our proof we use statement (\ref{eq:conv_mv})
that holds, in particular, for
    $$\varepsilon = \varepsilon^\prime = \frac{1}{4}\min\{ M(u)-M(v_0) \mid u\in B_{v_0}(3r),~ M(u)-M(v_0)>0 \}.$$
It means that, with probability $1$, there exists
$N=N(\varepsilon^\prime)$ such that for any $n>N$ and all $u\in
B_{v_0}(3r)$, we have $|M_n(u)-M(u)| < \varepsilon^\prime$.
Moreover, since $\ME(\xi_1) \subseteq B_{v_0}(3r)$, we can assert
the same for any $v\in\ME(\xi_1)$; namely, $|M_n(v)-M(v)| <
\varepsilon^\prime$. Together with the fact that $M(u)-M(v_0)>0$,
the obtained inequalities imply that, with probability $1$, there
exists number $N=N(\varepsilon^\prime)$ such that for any $n>N$ and
all $u\in B_{v_0}(3r)\setminus \ME(\xi_1)$,
    $$M_n(v_0) < M(v_0)+\varepsilon' \le M(v_0)+\frac{1}{4}( M(u)-M(v_0)) $$
    $$M(u)-\frac{1}{4}( M(u)-M(v_0)) \le M(u)-\varepsilon' < M_n(u),$$
and, hence,
    $$M_n(u)-M_n(v_0) \ge M(u)-\frac{1}{4}( M(u)-M(v_0)) - M(v_0)-\frac{1}{4}( M(u)-M(v_0)) = $$
    $$=\frac{1}{2}( M(u)-M(v_0)) \geq 2\varepsilon', ~ \mbox{i.e.,}$$
$$\BP \rbb{\exists N=N(\varepsilon),~~ \forall n>N,~~ \forall u\in B_{v_0}(3r) \setminus
\ME(\xi_1):~~ M_n(u)-M_n(v_0)\ge 2\varepsilon'} = 1.$$ Therefore,
(\ref{eq:main_eq_singleton_b}) holds for any $\delta \le
2\varepsilon'$. Choosing $\delta =
\min(\frac{1}{4}M(v_0),2\varepsilon')$ finishes the proof.

\end{proof}

\begin{corollary}[Inclusion Lemma] \label{le:slln_subset}
Let $\Gamma$ be a locally-finite connected graph,
$\{\xi_i\}_{i=1}^\infty$ a sequence of i.i.d. random
$\Gamma$-elements, and $\mu=\mu_{\xi_1}$.
Suppose that the weight function $M_{\xi}(\cdot)$ is totally defined. Then
    $$\BP\rb{\limsup_{n\rightarrow \infty}\MS(\xi_1,\ldots,\xi_n) \subseteq \ME (\xi_1) } = 1.$$
\end{corollary}

\begin{proof}
Lemma \ref{pr:separation_M} implies that
    $\BP\rbb{ u \notin \limsup \MS_n, \mbox{ for every } u\in V(\Gamma)\setminus\ME(\xi_1) }= 1.$
\end{proof}


\begin{theorem}\label{th:slln1}{\em(SLLN for graph-valued random elements with a singleton mean-set.)}
Let $\Gamma$ be a locally-finite connected graph and
$\{\xi_i\}_{i=1}^\infty$ a sequence of i.i.d. random
$\Gamma$-elements. If the weight function $M_{\xi_1}(\cdot)$ is totally defined and
    $\ME(\xi_1) = \{v\}$
for some $v\in V(\Gamma)$, then
    $$\lim_{n\rightarrow\infty} \MS(\xi_1,\ldots,\xi_n) = \ME(\xi_1)$$
almost surely (with probability one).
\end{theorem}

\begin{proof}
$\BP\rbb{\exists N \mbox{ s.t. }\forall n>N,~ M_n(v) < \inf_{u\in V(\Gamma)\setminus\{v\}} M_n(u) }  = 1,$
by Lemma \ref{pr:separation_M},
and, hence, $\BP\rbb{\exists N \mbox{ s.t. }\forall n>N,~ \MS(\xi_1,\ldots,\xi_n)=\{v\}}  = 1.$
\end{proof}

\subsection{Case of multi-vertex mean-sets}
\label{slln2or3vert} In this section we investigate a multi-vertex
mean-set case and conditions under which the strong law of large
numbers holds for such set. We reduce this problem to the question
of recurrence of a certain subset in $\MZ^n$ relative to a random
walk on this integer lattice. If $2 \leq |\ME(\xi)| \leq 3$, no
restrictive assumptions are required; we formulate and prove the law
for these special instances separately. The case $|\ME(\xi)| > 3$
requires more technical assumptions, and, thus, more work to handle
it.

\subsubsection{Preliminaries}

Assume $\ME(\xi_1) = \{v_1,v_2, \ldots, v_k\}$.
Our goal is to find conditions that would guarantee the
inclusion $\ME(\xi_1)\subseteq \limsup_{n\rightarrow \infty}\MS_n$
or, without loss of generality, conditions for $v_1 \in \limsup_{n\rightarrow \infty}\MS_n$.

By Lemma \ref{pr:separation_M}, it follows that, with probability
one, for a sequence of random $\Gamma$-elements $\{\xi_i\}_{i=1}^\infty$, there
exists a number $N$ such that for any $n>N$ we have
    $$\max\{M_n(v_1),M_n(v_2),\ldots,M_n(v_k)\} < \inf_{u\in\Gamma \setminus
\{v_1,v_2, \ldots, v_k\}}M_n(u).$$ Hence, for any $n>N$, $v_1 \in
\MS_n$ if and only if $M_n(v_1) \leq M_n(v_i)$ for every
$i=2,\ldots,k$. Thus, to achieve our goal, we need to show that
the system of inequalities
$$
\left\{
\begin{array}{l}
M_n(v_2) - M_n(v_1)\geq 0,\\
\ldots\\
M_n(v_k) - M_n(v_1) \geq 0,\\
\end{array}
\right.
$$
is satisfied for infinitely many $n\in \MN$.

For $i=1,\ldots,k-1$ and $n \in \MN$, define
    $$R_i(n) = n\rb{M_n(v_{i+1}) - M_n(v_1)} = \sum_{s\in \Gamma} \rb{d^2(v_{i+1},s)-d^2(v_1,s)} \cdot |\{i\mid \xi_i=s, 1\leq i\leq n\}|$$
and observe that
\begin{equation}\label{eq:step_R_i}
    R_i(n+1) - R_i(n) = \sum_{s\in \Gamma} ~~ [d^2(v_{i+1},s)-d^2(v_1,s)] ~~ \IND_{\{ \xi_{n+1}=s\}}.
\end{equation}
i.e., every $R_i(n)$ represents a random walk on $\MZ$ starting at $0$.
Consider a random walk $\ovR$, associated with $v_1$, in
$\MZ^{k-1}$, starting at the origin $(0,\ldots,0)$ with the position
of the walk after $n$ steps given by
    $$\ovR(n)= (R_1(n),R_2(n),\ldots,R_{k-1}(n)).$$
An increment step for $\ovR$ is defined by a vector $\ovzeta(s) =(\zeta_1(s), \ldots,\zeta_{k-1}(s)), ~ s\in V(\Gamma),$
with probability $\mu(s)$, where
    $$\zeta_i(s)=d^2(v_{i+1},s)-d^2(v_1,s).$$
The following lemma shows the significance of this random walk.

\begin{lemma}\label{le:ovR_recurence}
In the notation of this section, $v_1 \in \limsup_{n\rightarrow
\infty}\MS_n$ if and only if the random walk $\ovR$ visits the set
$\MZ_+^{k-1} = \{(a_1,\ldots,a_{k-1}) \mid a_i\ge 0\}$ infinitely
often. Therefore,
    $$\BP(v_1\in \limsup_{n\rightarrow \infty}\MS_n) = \BP(\ovR(n)\in \MZ_+^{k-1}, \mbox{ i.o.}).$$
\end{lemma}

\begin{proof}
Follows from the discussion preceding the lemma.
\end{proof}

It is worth redefining $\ovR$ in the terms of transition probability
function, as in \cite{Spitzer:2001}. Let $\ovo\in \MZ^{k-1}$ be the
zero vector and $x_i = \zeta_i(s)$, $s \in V(\Gamma)$. For every
$\ovx =(x_1, \ldots, x_{k-1}) \in \MZ^{k-1}$, we define a function
$P(\ovo,\ovx)$ by
\begin{equation}\label{eq:transition_probability}
    P(\ovo,\ovx) = \mu\{s \mid x_i = d^2(v_{i+1},s)-d^2(v_1,s) \mbox{ for every } i=1,\ldots,k-1\}.
\end{equation}
It is trivial to check that this is, indeed, the transition
probability for $\ovR$. To continue further, we investigate some
properties of our random walk $\ovR$.
\begin{lemma}\label{le:mu_m_finite}
Let $\ovR$ be a random walk defined above. Then
    $$m_1 = \sum_{\ovx\in\MZ^{k-1}} \ovx P(\ovo,\ovx) = \ovo ~~\mbox{ and }~~ m_2 = \sum_{\ovx\in\MZ^{k-1}} |\ovx|^2 P(\ovo,\ovx) < \infty.$$
\end{lemma}

\begin{proof}
The first equality trivially holds. For the second one, we get

    $$\sum_{\ovx\in\MZ^{k-1}} |\ovx|^2 P(\ovo,\ovx) = \sum_{s\in V(\Gamma)}\sum_{i=1}^{k-1} \rbb{d^2(v_{i+1},s)-d^2(v_{1},s)}^2  \mu(s)$$
    $$ \le \sum_{i=1}^{k-1} d^2(v_1,v_{i+1}) \sum_{s\in V(\Gamma)}\rbb{d(v_1,s)+d(v_{i+1},s)}^2 \mu(s)$$
    $$\le \sum_{i=1}^{k-1} d^2(v_1,v_{i+1}) (4M(v_1)+4M(v_{i+1})) < \infty.$$
\end{proof}

Clearly, conditions under which this random walk is recurrent would
guarantee that $v_1 \subseteq \limsup_{n\rightarrow \infty}\MS_n$
(see \cite[page 30, Proposition 3.3]{Spitzer:2001}). A general (not
simple, not symmetric) one-dimensional random walk is recurrent if
its first moment is zero and its first absolute moment is finite
(see \cite{Spitzer:2001}, pg. 23). Sufficient conditions for the
recurrence of two-dimensional random walk involve the finiteness of
its second moment and can be found in \cite[page 83]{Spitzer:2001}.
The result stated there indicates that genuinely 2-dimensional
random walk is recurrent if its first moment is zero, and its second
moment is finite. Let us recall some important notions before we go
on.

Consider an arbitrary random walk $\ovR$ on $\MZ^n$ given by a
transition probability $P$, as in (\ref{eq:transition_probability}).
The support, $supp(P)$, of the probability measure $P$ is defined to
be the set $supp(P) = \{\ovv \in \MZ^n \mid P(\ovv)\ne 0\}$ of all
possible one-step increments of $\ovR$. Further, with $\ovR$, one
can associate an abelian subgroup $A_\ovR$ of $\MZ^n$ generated by
the vectors in $supp(P)$. It is well-known in group theory that any
subgroup $A_\ovR$ of $\MZ^n$ is isomorphic to $\MZ^k$, where $k\le
n$ (the reader can also check \cite [Proposition7.1 on
pg.65]{Spitzer:2001} for details), in which case we write
$\dim(A_\ovR) = k$ and say that $\ovR$ is {\em genuinely
$k$-dimensional}.
Let us stress that we speak of an $n$-dimensional random walk on
$\MZ^n$ when $P(0,\ovx)$ is defined for all $\ovx$ in $\MZ^n$; this
walk is genuinely $n$-dimensional if $\dim(A_\ovR) = n$. We say that
$\ovR$ is {\em aperiodic} if $A_\ovR = \MZ^n$. Observe that
genuinely $n$-dimensional random walk does not have to be aperiodic.
A standard simple random walk, which we denote by $S = S(n)$, is an
example of an aperiodic random walk on $\MZ^n$. It will be
convenient to define a vector space $V_\ovR \subset \MR^n$ spanned
by the vectors in $supp(P)$. It is easy to see that the genuine
dimension of $\ovR$ is equal to the dimension of $V_\ovR$. We shall
need another notion for our developments. Assume that $D$ is an
$k\times n$ matrix (not necessarily integer valued) which maps
$A_{\ovR}$ onto $\MZ^k$. Then $D$ naturally induces a random walk
$\ovR^D$ on $\MZ^k$ with transition probability $P^D$ given by
$P^D(\ovu) = P( \ovv \in \MZ^n \mid D(\ovv)=\ovu )$ for every $\ovu
\in \MZ^k$.

\subsubsection{Strong law of large numbers for two or three vertices mean-sets}
Now, we can easily prove our strong law of large numbers for
mean-sets with two or three elements.

\begin{theorem}[SLLN for graph random elements with two or three point mean-set]\label{th:slln3}
Let $\Gamma$ be a locally-finite connected graph and
$\{\xi_i\}_{i=1}^\infty$ be a sequence of i.i.d.
random $\Gamma$-elements. If the weight
function $M_{\xi_1}(\cdot)$ is totally defined and $2 \leq
|\ME(\xi)| \leq 3$, then
    $$\limsup_{n\rightarrow \infty} \MS_n = \ME(\xi_1)$$
holds with probability one.
\end{theorem}

\begin{proof}
Assume that $\ME(\xi) = \{v_1,v_2\}$. Then the
random walk $\ovR$ is one-dimensional.
It is recurrent if
    $$\sum_{s \in \Gamma} | \zeta_1(s)|\mu(s) < \infty \mbox{ and } \sum_{s \in \Gamma} \zeta_1(s)\mu(s)= 0$$
(see \cite{Spitzer:2001}, pg. 23).
The equality $M(v_1) = M(v_2)$ implies the
second conditions and
    $$\sum_{s \in \Gamma} | \zeta_1(s)|\mu(s) = \sum_{s \in \Gamma} | d^2(v_2,s)-d^2(v_1,s) | \mu(s)$$
    $$\leq \sum_{s \in \Gamma} ( d^2(v_2,s) + d^2(v_1,s) ) \mu(s) = M(v_1) + M(v_2) < \infty$$
implies the first condition.
Hence, $\ovR$ is recurrent, and takes on
positive and negative values infinitely often.
We conclude that almost always
$\displaystyle \limsup_{n\rightarrow \infty} \MS_n = \{v_1,v_2\} = \ME \xi$.

Assume that $\ME(\xi) = \{v_1,v_2,v_3\}$.
Then the random walk $\ovR$ can be genuinely $0$, $1$, or $2$-dimensional.
The first case is trivial, the second can be considered as the case when $|\ME(\xi)|=2$.
So, assume $\ovR$ is genuinely $2$-dimensional.
By Lemma
\ref{le:mu_m_finite}, the first moment of $\ovR$ is $(0,0)$ and  the
second moment is finite. Now, it follows from \cite[Theorem
8.1]{Spitzer:2001} that $\ovR$ is recurrent. In particular,
$\MZ_+^{k-1}$ is visited infinitely often
with probability $1$.

In both cases, it follows from Lemma \ref{le:ovR_recurence}
that $\BP(v_1 \in \limsup_{n\rightarrow \infty} \MS_n) = 1$.
Hence the result.
\end{proof}

Recall that a subset of $\MZ^n$ is called {\em recurrent} if it is visited
by a given random walk infinitely often with probability one, and it
is transient otherwise.
A criterion for recurrence of a set for a
simple random walk was obtained in \cite{ItoMcKean:1960} for $n=3$
(it can also be found in \cite[Theorem 26.1]{Spitzer:2001}). It
turns out that the criterion does not depend on a random walk in
question. This is the subject of the extension of the Wiener's test,
proved in \cite{Uchiyama:1998}, that we state below. This invariance
principle is one of the main tools we use in our investigation of
the recurrence properties of the positive octant in $\MZ^n$ for
$\ovR$.

\medskip{\noindent \bf Theorem.}  (Extension of Wiener's test,
\cite{Uchiyama:1998}) {\em Let $n\ge 3$. An infinite subset
$A$ of $\MZ^n$ is either recurrent for each aperiodic random walk
$\ovR$ on $\MZ^n$ with mean zero and a finite variance, or transient
for each of such random walks.}

For a positive constant $\alpha \in\MR$ and a positive integer $m\le
n$ define a subset of $\MR^n$
    $$Cone_\alpha^m = \cb{(x_1,\ldots,x_n) \in \MR^n \mid x_1=0,\ldots,x_{n-m}=0,~\sqrt{x_{n-m+1}^2+\ldots+x_{n-1}^2}\le \alpha x_n}$$
called an $m$-dimensional cone in $\MR^n$. If $m=n$, then we omit
the superscript in $Cone_\alpha^m$. For an $n\times n$ matrix $D$
and a set $A\subseteq \MR^n$, define a set $A^D = \{D\cdot\ovv \mid
\ovv\in A\}$, which is a linear transformation of $A$. If $D$ is an
orthogonal matrix, then the set $(Cone_\alpha)^D$ is called a {\em
rotated cone}. As in \cite{ItoMcKean:1960}, for any non-decreasing
function $i:\MN \rightarrow \MR_+$ define a set
    $$Thorn_i = \{\ovv\in\MZ^n \mid \sqrt{v_1^2 +\ldots + v_{n-1}^2}\le i(v_n)\}.$$
Observe that $Cone_\alpha \cap \MZ^n = Thorn_i$ where $i(t) = \alpha
t$.

\begin{theorem}\label{th:cone_visit}
For any $\alpha>0$ and any orthogonal matrix $D$,
    $$\BP\rb{S(n)\in (Cone_\alpha)^D, \mbox{ i.o.}} =1,$$
i.e.,  the probability that the simple random walk on $\MZ^n$ visits
$(Cone_\alpha)^D$ infinitely often is $1$.
\end{theorem}

\begin{proof}
Direct consequence of (6.1) and (4.3) in \cite{ItoMcKean:1960},
where the criterion for recurrence of $Thorn_i$ is given.
\end{proof}

Next two lemmas are obvious

\begin{lemma}\label{le:D_rotates_cones}
Assume that a set $A\subseteq \MR^n$ contains a rotated cone. Then
for any invertible $n\times n$ matrix $D$, the set $A^D$ contains a
rotated cone.
\end{lemma}

\begin{lemma}
If $S_1\subseteq S_2 \subseteq \MR^n$ and $S_1$ is visited by the
simple random walk infinitely often with probability $1$ then $S_2$
is visited by the simple random walk infinitely often with
probability $1$.
\end{lemma}

Now, we return to our strong law of large numbers for multi-vertex
mean-sets. Assume that $\ME \xi = \{v_1,\ldots, v_k\}$, where $k\ge
4$. Let $\ovR^i$ be a random walk on $\MZ^{k-1}$, associated with
$v_i$, where $i=1,\ldots,k$ (in our notation, $\ovR = \ovR^1$ ).
This is a $(k-1)$-dimensional random walk which, in general, is not
aperiodic. In fact, $\ovR^i$ is not even genuinely
$(k-1)$-dimensional. Fortunately, it turns out that it does not
matter to what vertex $v_i$ we associate our random walk, since the
choice of the vertex does not affect the dimension of the
corresponding walk, as the following lemma shows.

\begin{lemma}\label{le:equal_rw_dimensions}
Let $\mu$ be a probability measure on a locally finite graph
$\Gamma$ such that $\ME \mu = \{v_1,\ldots,v_k\}$, where $k\ge 2$.
Then the random walks $\ovR^1,\ldots,\ovR^{k}$, associated with
vertices $v_1,\ldots,v_k$ respectively, all have the same genuine
dimension.
\end{lemma}

\begin{proof}
We prove that random walks $\ovR^1$ and $\ovR^{2}$ have the same
genuine dimension. Recall that the subgroup $A_{\ovR^1}$ is
generated by the set of vectors $\ovv^1\in\MZ^{k-1}$ such that for
some $s\in supp(\mu)$, $\ovv^1 = \ovv^1(s) =
(d^2(v_2,s)-d^2(v_1,s),d^2(v_3,s)-d^2(v_1,s),\ldots,d^2(v_k,s)-d^2(v_1,s))$
and the subgroup $A_{\ovR^2}$ is generated by the set of vectors
$\ovv^2\in\MZ^{k-1}$ such that for some $s\in supp(\mu)$, $\ovv^2 =
\ovv^2(s) =
(d^2(v_1,s)-d^2(v_2,s),d^2(v_3,s)-d^2(v_2,s),\ldots,d^2(v_k,s)-d^2(v_1,s))$.
Observe that for every $s\in supp(\mu)$ the equality $\ovv^2(s)
= D \cdot \ovv^1(s)$ holds, where $D$ is a $(k-1)\times (k-1)$ matrix
$$
D = \left(
\begin{array}{ccccc}
-1 & 0 & 0 & 0 & \ldots \\
-1 & 1 & 0 & 0 & \ldots\\
-1 & 0 & 1 & 0 & \ldots\\
-1 & 0 & 0 & 1 & \ldots\\
\ldots \\
\end{array}
\right)
$$
Therefore, $A_{\ovR^2} = (A_{\ovR^1})^D$. Since the matrix $D$ is
invertible it follows that $A_{\ovR^1}$ and $A_{\ovR^2}$ have the
same dimension.
\end{proof}

\begin{theorem}\label{th:slln4}
Let $\Gamma$ be a locally-finite connected graph and
$\{\xi_i\}_{i=1}^\infty$ a sequence of i.i.d.
random $\Gamma$-elements. Assume that the
weight function $M$ is totally defined and $\ME(\xi) =
\{v_1,\ldots,v_k\}$, where $k\ge 4$. If the random walk $\ovR^1$
associated to $v_1$ is genuinely $(k-1)$-dimensional, then
    $$\limsup_{n\rightarrow \infty}\MS_n = \ME(\xi_1)$$
holds with probability $1$.
\end{theorem}

\begin{proof}
Since $\ovR^1$ is genuinely $(k-1)$-dimensional it follows that the
subgroup $A_{\ovR^1}$ is isomorphic to $\MZ^{k-1}$ and there exists
an invertible matrix $D$ that isomorphically maps $A_{\ovR^1}
\subseteq \MZ^{k-1}$ onto $\MZ^{k-1}$. Consider a set $\MR_+^{k-1} =
\{(x_1,\ldots,x_{k-1}) \mid x_i\ge 0\}$. Obviously,
$\BP\rb{\ovR^1\in\MZ_+^{k-1} \mbox{ i.o.}} =
\BP\rb{\ovR^1\in\MR_+^{k-1} \mbox{ i.o.}}$.

Let $(\ovR^1)^D$ be the random walk on $\MZ^{k-1}$ induced by $D$ by
application of $D$ to $\ovR^1$. The random walk $(\ovR^1)^D$ is
aperiodic since $D$ maps $A_{\ovR^1}$ onto $\MZ^{k-1}$ and, by
construction of $(\ovR^1)^D$,
    $$\BP\rb{\ovR^1\in\MR_+^{k-1} \mbox{i.o.}} = \BP\rb{(\ovR^1)^D\in(\MR_+^{k-1})^D \mbox{ i.o.}}.$$
Let $S$ be the simple random walk on $\MZ^{k-1}$. Since $(\ovR^1)^D$ and $S$
are both aperiodic random walks on $\MZ^{k-1}$, it follows from the
Invariance Principle (Extension of Wiener's test) that
    $$\BP\rb{(\ovR^1)^D\in(\MR_+^{k-1})^D \mbox{ i.o.}} = \BP\rb{S\in(\MR_+^{k-1})^D \mbox{ i.o.}}.$$
Clearly, the set
$\MR_+^{k-1}$ contains a rotated cone and, hence, by Lemma
\ref{le:D_rotates_cones}, its image under an invertible linear
transformation $D$ contains a rotated cone too. Now, by Theorem
\ref{th:cone_visit},
$\BP\rb{S\in(\MR_+^{k-1})^D \mbox{ i.o.}} = 1$.
Thus, $\BP\rb{\ovR^1\in \MZ_+^{k-1} \mbox{ i.o.}} = 1$ and by Lemma
\ref{le:ovR_recurence}
    $$\BP(v_1\in\limsup_{n\rightarrow \infty}\MS_n) = 1.$$
Finally, it follows from Lemma \ref{le:equal_rw_dimensions} that for
any $i=2,\ldots,k$ the random walk $\ovR^i$ is genuinely
$(k-1)$-dimensional. For any $i=2,\ldots,k$ we can use the same
argument as for $v_1$ to prove that $\BP(v_i\in\limsup_{n\rightarrow
\infty}\MS_n) = 1.$ Hence the result.
\end{proof}

\subsubsection{The case when random walk is not genuinely $(k-1)$-dimensional}

The case when $\ovR^1$ is not genuinely $(k-1)$-dimensional is more
complicated. To answer the question whether $v_1$ belongs to
$\limsup_{n\rightarrow \infty}\MS_n$ (namely, how often $v_1\in \MS_n$),
we need to analyze how the space $V_{\ovR^1}$
``sits'' in $\MR^{k-1}$. We know that the subgroup $A_{\ovR^1}
\subset \MZ^{k-1}$ is isomorphic to $\MZ^m$, where $m<k-1$ in the
case under consideration. Therefore, there exists a $m\times (k-1)$
matrix $D$ which maps the subgroup $A_{\ovR^1}$ onto $\MZ^m$ and
which is injective onto $A_{\ovR^1}$. Furthermore, the mapping $D$
maps the subspace $V_{\ovR^1}$ bijectively onto $\MR^{m}$. The
linear mapping $D$ induces an aperiodic random walk $(\ovR^1)^D$ on
$\MZ^m$ in a natural way and $\BP\rb{\ovR^1\in(\MR_+^{k-1}) \mbox{
i.o.}} = \BP\rb{\ovR^1\in(\MR_+^{k-1}\cap V_{\ovR^1}) \mbox{ i.o.}}
= \BP\rb{(\ovR^1)^D\in(\MR_+^{k-1}\cap V_{\ovR^1})^D \mbox{ i.o.}}$.
The main problem here is to understand the structure of the set
$(\MR_+^{k-1}\cap V_{\ovR^1})^D$ and, to be more precise, the
structure of the set $B_{\ovR^1} = \MR_+^{k-1}\cap V_{\ovR^1}$.
Clearly $B_{\ovR^1}$ is a monoid, i.e., contains the trivial element
and a sum of any two elements in $B_{\ovR^1}$ belongs to
$B_{\ovR^1}$. We can define dimension of $B_{\ovR^1}$ to be the
maximal number of linearly independent vectors in $B_{\ovR^1}$.

\begin{theorem} \label{th:dimension}
Suppose $A_{\ovR^1}\simeq \MZ^m$ and the set $B_{\ovR^1}$ has
dimension $m$. Then $\BP(v_i\in\limsup_{n\rightarrow \infty}\MS_n) =
1.$
\end{theorem}

\begin{proof}
Since $B_{\ovR^1}$ is a monoid of dimension $m$ it is not hard to
see that $B_{\ovR^1}$ contains an $m$-dimensional rotated cone.
Since $D$ is a linear isomorphism from $V_{\ovR^1}$ onto $\MR^m$ it
follows by Lemma \ref{le:D_rotates_cones} that $(B_{\ovR^1})^D$
contains an $m$-dimensional rotated cone in $\MR^m$. If $S$ is a
simple random walk in $\MZ^m$ then
    $\BP\rb{S\in(B_{\ovR^1})^D \mbox{ i.o.}} = 1$
and since $S$ and $(\ovR^1)^D$ are both aperiodic, by the extension
of Wiener's test (Invariance Principle), we see that
    $\BP\rb{(\ovR^1)^D\in(B_{\ovR^1})^D \mbox{ i.o.}} = 1$.
Hence, $\BP\rb{\ovR^1\in(\MR_+^{k-1}) \mbox{ i.o.}} = 1$ and by
Lemma \ref{le:ovR_recurence}, $\BP(v_i\in\limsup_{n\rightarrow
\infty}\MS_n) = 1.$

\end{proof}

Below we investigate under what conditions the subgroup $A_{\ovR^1}$
and the set $B_{\ovR^1}$ have the same dimension $m$.

\begin{lemma}\label{le:posit_inters}
Assume that $A_{\ovR^1}$ contains a positive vector. Then
$A_{\ovR^1}$ and the set $\MR_+^{k-1} \cap V_{\ovR^1}$ have the same
dimension.
\end{lemma}

\begin{proof}
Straightforward.
\end{proof}

\begin{lemma}\label{le:posit_mu}
Assume that $\mu(v_1) \ne 0$. Then $A_{\ovR^1}$ and the set
$\MR_+^{k-1} \cap V_{\ovR^1}$ have the same dimension.
\end{lemma}

\begin{proof}
Observe that if $\mu(v_1) \ne 0$ then $A_{\ovR^1}$ contains the
vector $(d^2(v_2,v_1),\ldots,d^2(v_k,v_1))$ which has all positive
coordinates. Therefore, by Lemma \ref{le:posit_inters} the set
$A_{\ovR^1}$ and $\MR_+^{k-1} \cap V_{\ovR^1}$ have the same
dimension.
\end{proof}

\begin{corollary}\label{th:slln5}
Let $\Gamma$ be a locally-finite connected graph and
$\{\xi_i\}_{i=1}^\infty$ be a sequence of i.i.d.
random $\Gamma$-elements. Assume that the
weight function $M_{\xi_1}(\cdot)$ is totally defined and $\ME(\xi)
= \{v_1,\ldots,v_k\}$, where $k\ge 4$. If $\ME (\xi_1) \subseteq
supp(\mu)$ then $\limsup_{n\rightarrow \infty}\MS_n = \ME(\xi_1)$
holds with probability one.
\end{corollary}

\begin{proof}
Follows from Lemma \ref{le:posit_mu}, \ref{le:posit_inters}, and
\ref{th:dimension}.
\end{proof}

\section{Concentration of measure inequalities}
\label{se:Chebyshev}

Concentration inequalities are upper bounds on the rate of convergence (in probability) of
sample (empirical)
means to their ensemble counterparts (actual means).
Chebyshev inequality and Chernoff-Hoeffding exponential bounds
are classical examples of such inequalities in probability theory.
In this section, we prove analogues of the classical Chebyshev's
inequality and Chernoff-Hoeffding like bounds - the concentration of measure
inequalities for a graph-
(group-)valued random elements.

\subsection{Chebyshev's inequality for graphs/groups}

The classical Chebyshev's inequality asserts that if $\xi$ is a
random variable with $\ME(\xi ^2)<\infty$, then for any $\varepsilon
> 0$, we have
\begin{equation}\label{eq:Cheb_ineq1}
\BP(|\xi - \ME(\xi)|\ge \varepsilon) \le \frac{\sigma^2}{\varepsilon
^2},
\end{equation}
where $\sigma^2 = Var(\xi)$, see \cite{Billingsley}.

Chebyshev discovered it when he was trying to prove the law of large
numbers, and the inequality is widely used ever since. Chebyshev's inequality is
a result concerning the concentration
of measure, giving a quantitative description of this concentration.
Indeed, it provides a bound on the probability that a value of a
random variable $\xi$ with finite mean and variance will differ from
the mean by more than a fixed number $\varepsilon$. In other words,
we have a crude estimate for concentration of probabilities around
the expectation, and this estimate has a big theoretical
significance.

The inequality
(\ref{eq:Cheb_ineq1}) applied to the sample mean random variable
$\ovX = \frac{S_n}{n}$, where $S_n = \xi_1 + \ldots + \xi_n$,
$\ME(\xi_i)=m, Var(\xi_i)= \sigma ^2, i=1,\ldots,n$ results in
\begin{equation}\label{eq:Cheb_ineq2}
\BP(|\ovX - m|\ge \varepsilon) \le \frac{\sigma^2}{n \varepsilon ^2}
\end{equation}
The goal is to prove a similar inequality for a graph-valued
random element $\xi$.

\begin{lemma}\label{le:_mu_M_cond}
Let $\mu$ be a distribution on a locally finite graph $\Gamma$ such
that $M\equiv M^{(2)}$ is defined. If for some $r\in \MN$ and $v_0
\in V(\Gamma)$ the inequality
\begin{equation}\label{eq:finiteE_cond}
    \sum_{s\in V(\Gamma) \setminus B_{v_0}(r/2)} d(v_0,s) \mu(s) - \frac{r}{2} \mu(v_0) < 0
\end{equation}
holds, then $M(u) > M(v_0)$ for any $u \in V(\Gamma) \setminus
B_{v_0}(r)$.
\end{lemma}

\begin{proof}
Indeed, pick any $u \in V(\Gamma) \setminus B_{v_0}(r)$ and put $d =
d(v_0,u)$. Then
    $$M(u)-M(v_0) = \sum_{s\in V(\Gamma)} (d^2(u,s)-d^2(v_0,s))\mu(s)$$
    $$\ge d^2\mu(v_0) - \sum_{d(v_0,s)>d(u,s)} (d^2(v_0,s)-d^2(u,s))\mu(s)$$
    $$\ge d^2\mu(v_0) - 2d\sum_{d(v_0,s)>d(u,s)} d(v_0,s)\mu(s)
    \ge d^2\mu(v_0) - 2d\sum_{s \in V(\Gamma) \setminus B_{v_0}(r/2)} d(v_0,s)\mu(s).$$
Since $d>r$ it follows that the last sum is positive. Thus
$M(u)>M(v_0)$ as required.

\end{proof}

\begin{theorem}\label{th:Chebyshev}
Let $\Gamma$ be a locally-finite connected graph and
$\{\xi_i\}_{i=1}^\infty$ a sequence of i.i.d. random
$\Gamma$-elements. If the weight function
$M_{\xi_1}$ is totally defined and
    $\ME(\xi_1) = \{v\}$
for some $v\in V(\Gamma)$ then there exists a constant $C =
C(\Gamma,\xi_1)>0$ such that
\begin{equation}\label{eq:Chebyshev}
    \BP\rbb{ \MS(\xi_1,\ldots, \xi_n) \ne \{v\} } \le \frac{C}{n}.
\end{equation}
\end{theorem}

\begin{proof}
It follows from the definition of the sample mean-set that
    $$\{ \MS_n \ne \{v\} \} = \{\exists u\in V(\Gamma) \setminus \{v\},~~ M_n(u) \le M_n(v) \}.$$
Hence, it is sufficient to prove that $\BP \biggl(\exists u\in
V(\Gamma) \setminus \{v\},~~ M_n(u)\le M_n(v)\biggr)\le
\frac{C}{n}$, for some constant $C$. We do it in two stages. We show
that for some $v_0\in V(\Gamma)$ and constants $r \in \MN$, $C_1,C_2
\in \MR$ such that $v\in B_{v_0}(r)$ and inequalities
\begin{equation}\label{eq:Chbysh_part1}
    \BP \biggl(\exists u\in B_{v_0}(r) \setminus \{v\},~~ M_n(u)\le M_n(v)\biggr)\le \frac{C_1}{n}
\end{equation}
and
\begin{equation}\label{eq:Chbysh_part2}
    \BP\biggl(\exists u\in V(\Gamma) \setminus B_{v_0}(r),~~ M_n(u)\le M_n(v_0)\biggr)\le \frac{C_2}{n}
\end{equation}
hold. Clearly, for any $u,v_0,v\in V(\Gamma)$ if $M_n(u)\le M_n(v)$
then either $M_n(u) \le M_n(v_0)$ or $M_n(v_0)\le M_n(v)$. It is not
hard to see that if we find $C_1$ and $C_2$ satisfying
(\ref{eq:Chbysh_part1}) and (\ref{eq:Chbysh_part2}) respectively,
then (\ref{eq:Chebyshev}) holds for $C = C_1+C_2$ and the theorem is
proved.

First we argue (\ref{eq:Chbysh_part2}). Choose any $v_0 \in
V(\Gamma)$ such that $\mu(v_0)>0$ and $r\in \MN$ such that the
inequality (\ref{eq:finiteE_cond}) holds. We can choose such $r$
since $M^{(1)}(v_0)$ is defined. Observe that the left hand side of
the inequality above is the expectation of a random variable $\eta:
V\rightarrow \MR$ defined as $\eta(s) = d(v_0,s) \IND_{V(\Gamma)
\setminus B_{v_0}(r/2)}(s) - \frac{r}{2}\IND_{v_0}(s)$, $s \in
V(\Gamma)$, where $\IND_{\cdot}(s)$ is an indicator function. Since
by our assumption $M\equiv M^{(2)}$ is defined, it follows that
$\sigma^2(\eta)$ is defined, and, applying Lemma \ref{le:_mu_M_cond}
and the Chebyshev inequality with $\varepsilon=|\ME\eta|/2$, we
obtain
    $$\BP\biggl(\exists u\in V(\Gamma) \setminus B_{v_0}(r),~~ M_n(u)\le M_n(v_0)\biggr)$$
    $$\le\BP \biggl(\modbb{\sum_{s\in V(\Gamma) \setminus B_{v_0}(r/2)} d(v_0,s) \mu_n(s) - \frac{r}{2} \mu_n(v_0)-\ME \eta}\ge |\ME\eta|/2 \biggr)
    \le \frac{4\sigma^2(\eta)}{n|\ME\eta|^2}.$$
Hence, inequality (\ref{eq:Chbysh_part2}) holds for $C_2=C_2(r, v_0,
\mu) = \frac{4\sigma^2(\eta)}{|\ME\eta|^2}$.
To prove (\ref{eq:Chbysh_part1}) we notice that for any $u\in
V(\Gamma) \setminus \{v\}$,
    $\displaystyle M(u) - M(v) = \sum_{s\in V(\Gamma)} (d(u,s) - d(v,s))(d(u,s) + d(v,s)) \mu(s),$
i.e., $M(u)-M(v)$ is the expectation of a random variable $\tau:
V\rightarrow \MR$ defined as
    $$\tau_{u,v}(s)=(d(u,s) - d(v,s))(d(u,s) + d(v,s)), ~ s\in V(\Gamma).$$
Furthermore, since $M_{\xi_1}(\cdot)$ is defined and  for every
$s\in V(\Gamma)$, $d(u,s) - d(v,s) \le d(v,u)$, it is easy to see
that $\sigma^2(\tau_{u,v}(s))<\infty$. Thus, by the Chebyshev
inequality for the sample average of $\tau_{u,v}(s)$,
    $$\BP \biggl(|M_n(u)-M_n(v)-(M(u)-M(v))|\ge \varepsilon \biggr) \le \frac{\sigma^2(\tau_{u,v}(s))}{n\varepsilon^2}.$$
holds. Now, if $0<\varepsilon < M(u)-M(v)$, then
    $$\BP \biggl(M_n(u)<M_n(v) \biggr) \le \BP \biggl(|M_n(u)-M_n(v)-(M(u)-M(v))|\ge \varepsilon \biggr).$$
Finally, we choose $\varepsilon$ to be $\frac{1}{2} \inf\{M(u)-M(v)
\mid u\in B_{v_0}(r) \setminus \{v\} \}$ and using
$\sigma$-additivity of measure we see that inequality
(\ref{eq:Chbysh_part1}) holds for the constant $C_1 =
\varepsilon^{-2}\sum_{u\in B_{v_0}(r)} \sigma^2(\tau_{u,v}(s))$.

\end{proof}

\noindent In fact, one can easily generalize the previous theorem to
the following statement.

\begin{theorem}
Let $\Gamma$ be a locally-finite connected graph and
$\{\xi_i\}_{i=1}^\infty$ a sequence of i.i.d. random
$\Gamma$-elements. If the weight function
$M_{\xi_1}$ is totally defined then there exists a constant
$C = C(\Gamma,\xi_1)>0$ such that
\begin{equation}\label{eq:Chebyshev2}
    \BP\rbb{ \MS(\xi_1,\ldots, \xi_n) \not\subseteq \ME(\xi) } \le \frac{C}{n}.
\end{equation}
\end{theorem}

\subsection{Chernoff-Hoeffding like bound for graphs/groups}

Let $x_i$ be independent random variables. Assume that each $x_i$ is
almost surely bounded, i.e., assume that for every $i\in\MN$ there
exists $a_i,b_i\in \MR$ such that $\BP(x_i-\ME x_i \in [a_i,b_i])
=1$. Then for $S_n = \sum_{i=1}^n x_i$ and for any $\varepsilon>0$
we have the inequality (called the Hoeffding's inequality)
    $$\BP(|S_n-\ME S_n| \ge n\varepsilon) \le 2 \exp\rb{-\frac{2n^2\varepsilon^2}{\sum_{i=1}^n (b_i-a_i)^2}}.$$
If $x_i$ are identically distributed then we get the inequality
    $$\BP\rb{\modb{\frac{1}{n}(x_1+\ldots+x_n) - \ME x_1} \ge \varepsilon} \le 2 \exp\rb{-\frac{2\varepsilon^2}{(b-a)^2}n}.$$
Techniques of the previous section can be used to find a similar
bound on $\BP\rbb{ \MS(\xi_1,\ldots, \xi_n) \not\subseteq \ME(\xi)
}$ for a sequence of iid graph-valued $\xi_i$ satisfying some
additional assumptions.

\begin{theorem}
Let $\Gamma$ be a locally-finite connected graph and
$\{\xi_i\}_{i=1}^\infty$ a sequence of i.i.d. random
$\Gamma$-elements. If the weight function
$M_{\xi_1}(\cdot)$ is totally defined and $\mu_{\xi_1}$ has finite
support then for some constant $C>0$
\begin{equation}\label{eq:Hoeffding2}
    \BP\rbb{ \MS(\xi_1,\ldots, \xi_n) \not\subseteq \ME(\xi) } \le O(e^{-Cn}).
\end{equation}
\end{theorem}

\begin{proof}
Proof is similar to the proof of Theorem \ref{th:Chebyshev}. We find
$v_0\in V(\Gamma)$, $r\in \MN$, and constants $C_1,C_2>0$ such that
inequalities
\begin{equation}\label{eq:Hoeffding_part1}
    \BP \biggl(\exists u\in B_{v_0}(r) \setminus \{v\},~~ M_n(u)\le M_n(v)\biggr)\le O(e^{-C_1n})
\end{equation}
and
\begin{equation}\label{eq:Hoeffding_part2}
    \BP\biggl(\exists u\in V(\Gamma) \setminus B_{v_0}(r),~~ M_n(u)\le M_n(v_0)\biggr)\le O(e^{-C_2n})
\end{equation}
hold.

Choose $v_0\in V(\Gamma)$ and $r\in \MN$ exactly the same way as in
Theorem \ref{th:Chebyshev}. Note that a random variable $\eta(s) =
d(v_0,s) \IND_{V(\Gamma) \setminus B_{v_0}(r/2)}(s) -
\frac{r}{2}\IND_{v_0}(s)$ (where $s \in V(\Gamma)$) is almost surely
bounded. Choose a lower and an upper bounds for $\eta$ and denote
them by $a$ and $b$ respectively. Now, applying Hoeffding's
inequality to $\eta$ with $\varepsilon = |\ME \eta|/2$ we obtain
    $$\BP\biggl(\exists u\in V(\Gamma) \setminus B_{v_0}(r),~~ M_n(u)\le M_n(v_0)\biggr)$$
    $$\le\BP \biggl(\modbb{\sum_{s\in V(\Gamma) \setminus B_{v_0}(r/2)} d(v_0,s) \mu_n(s) - \frac{r}{2} \mu_n(v_0)-\ME \eta}\ge |\ME\eta|/2 \biggr)
    \le 2 \exp \rb{-\frac{|\ME \eta|^2}{2(b-a)^2}n}.$$
Therefore, (\ref{eq:Hoeffding_part2}) holds for $C_2 = \frac{|\ME
\eta|^2}{2(b-a)^2}$.

To prove (\ref{eq:Hoeffding_part1}) we notice that for any
$v\in\ME(\xi)$ and $u\in V(\Gamma) \setminus \ME(\xi)$ we have
    $\displaystyle M(u) - M(v) = \sum_{s\in V(\Gamma)} (d(u,s) - d(v,s))(d(u,s) + d(v,s)) \mu(s),$
i.e., $M(u)-M(v)$ is the expectation of a random variable
$\tau_{u,v}: V\rightarrow \MR$ defined as
    $$\tau_{u,v}(s)=(d(u,s) - d(v,s))(d(u,s) + d(v,s)), ~ s\in V(\Gamma).$$
Furthermore, since $\xi_1$ has finite support it follows that the
random variable $\tau_{u,v}(s)$ almost surely belongs to
$[a_{u,v},b_{u,v}]$. Thus, by the Hoeffding's inequality for the
sample average of $\tau_{u,v}(s)$,
    $$\BP \biggl(|M_n(u)-M_n(v)-(M(u)-M(v))|\ge \varepsilon \biggr) \le 2\exp\rb{-\frac{2\varepsilon^2}{(b_{u,v}-a_{u,v})^2}n}.$$
holds. Now, if $0<\varepsilon < M(u)-M(v)$, then
    $$\BP \biggl(M_n(u)<M_n(v) \biggr) \le \BP \biggl(|M_n(u)-M_n(v)-(M(u)-M(v))|\ge \varepsilon \biggr).$$
Choose $\varepsilon$ to be $\frac{1}{2} \inf\{M(u)-M(v) \mid v\in
\ME(\xi),~ u\in B_{v_0}(r) \setminus \ME(\xi) \}$ and $\delta =
\max\{b_{u,v}-a_{u,v} \mid v\in \ME(\xi),~ u\in B_{v_0}(r) \setminus
\ME(\xi)\}$. Finally, using $\sigma$-additivity of measure we see
that inequality (\ref{eq:Hoeffding_part1}) holds for the constant
$C_1 = \frac{2\varepsilon^2}{\delta^2}$.
\end{proof}

\section{Configurations of mean-sets with applications}\label{se:conf_centers}

In this section, we discuss several configurations of mean-sets
on graphs and, in particular, on trees and free groups.
First, we make a simple observation stated in the lemma below.

\begin{lemma}
Let $\Gamma$ be a connected graph. Then for any $v \in V(\Gamma)$
there exists a measure $\mu$ such that $\ME (\mu) = \{v\}$.
\end{lemma}

\begin{proof}
Indeed, the statement of the lemma holds for the distribution
defined by
$$
\mu(u) = \left\{
\begin{array}{ll}
1, & \mbox{if } u=v;\\
0, & \mbox{otherwise}.\\
\end{array}
\right.
$$
\end{proof}

\begin{figure}[h]
\centerline{\includegraphics[scale=1]{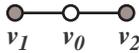} }
\caption{Impossible configuration of centers (gray vertices).}
\label{fi:graph2}
\end{figure}

On the other hand, it is easy to see that not any subset of
$V(\Gamma)$ can be realized as $\ME (\mu)$. For instance, consider a
graph as in Figure \ref{fi:graph2}.
Let $\mu_0=\mu(v_0)$, $\mu_1=\mu(v_1)$, $\mu_2=\mu(v_2)$,
$M_0=M(v_0)$, $M_1=M(v_1)$, $M_2=M(v_2)$ Then
    $M_1 = \mu_0+4\mu_2,$
    $M_0 = \mu_1+\mu_2,$
    $M_2 = 4\mu_1+\mu_0.$
Clearly, for no values of $\mu_0$, $\mu_1$, and $\mu_2$ both
inequalities $M_0>M_1$ and $M_0>M_2$ can hold simultaneously (since
we can not have $2M_0>M_1 + M_2$). Thus, $v_1$ and $v_2$ can not
comprise $\ME\mu$. In fact, a tree can have only a limited
configuration of centers as proved in Proposition
\ref{pr:tree_centers} below.

Let $\Gamma$ be a graph. We say that $v_0\in V(\Gamma)$ is a
 \textit{cut-point} if removing $v_0$ from $\Gamma$ results into a
disconnected graph.
The same definition holds for any metric space. It turns out that
existence of a \textit{cut-point} in $\Gamma$ affects configurations
of mean-sets. The following lemma provides a useful inequality that
holds for any metric space with a cut-point.

\begin{lemma}[Cut-point inequality] \label{le:cut_point_ineq1}
Let $(\Gamma,d)$ be a metric space and $v_0$ a cut point in
$\Gamma$. If $v_1$, $v_2$ belong to distinct connected components of
$\Gamma \setminus \{v_0\}$ then for any $s \in V(\Gamma)$ the
inequality
\begin{equation}\label{eq:cut_point_ineq}
\begin{array}{c}
d(v_0,v_2) \rb{d^2(v_1,s)-d^2(v_0,s)} + d(v_0,v_1) \rb{d^2(v_2,s)-d^2(v_0,s)}\ge C > 0 \\

\end{array}
\end{equation}
holds, where $C = C(v_0, v_1, v_2)= d(v_0,v_2)d(v_0,v_1)(d(v_0,v_1)+
d(v_0,v_2) )$.
\end{lemma}

\begin{proof}
Denote the left hand side of (\ref{eq:cut_point_ineq}) by $g(s)$.
There are $3$ cases to consider.

{\sc Case 1.} Assume that $s$ does not belong to the components of
$v_1$ and $v_2$. Then
    $$d(v_0,v_2) \rb{d^2(v_1,s)-d^2(v_0,s)} + d(v_0,v_1) \rb{d^2(v_2,s)-d^2(v_0,s)}$$
    $$= d(v_0,v_2) d(v_0,v_1) \rb{2d(v_0,s) + d(v_0,v_1)} + d(v_0,v_1) d(v_0,v_2) \rb{2d(v_0,s) + d(v_0,v_2)}$$
    $$= d(v_0,v_2) d(v_0,v_1) \rb{4d(v_0,s) + d(v_0,v_1) + d(v_0,v_2)}
    \ge d(v_0,v_2) d(v_0,v_1) \rb{d(v_0,v_1) + d(v_0,v_2)}$$
and hence (\ref{eq:cut_point_ineq}) holds.

{\sc Case 2.} Assume that $s$ belongs to the component of $v_1$.
Define
    $$x = x(s) = d(v_1,s) \mbox{ and } y = y(s) = d(v_0,s).$$
In this notation we get
    $$g(s) = g(x,y) = d(v_0,v_2) \rb{x^2-y^2} + d(v_0,v_1) \rb{2yd(v_0,v_2)+d^2(v_0,v_2)}.$$
Dividing by a positive value $d(v_0,v_2)$, we get
    $$g(s) >0 ~~\mbox{ if and only if }~~ \frac{g(x,y)}{d(v_0,v_2)} = x^2-y^2 + d(v_0,v_1) \rb{2y+d(v_0,v_2)} > 0.$$
Now, observe that the numbers $x$, $y$, and $d(v_0,v_1)$ satisfy
triangle inequalities
    $$
\left\{
\begin{array}{l}
x+y \ge d(v_0,v_1);\\
x+d(v_0,v_1) \ge y;\\
y+d(v_0,v_1) \ge x;\\
\end{array}
\right.
    $$
that bound the area visualized in Figure \ref{fi:triangle_area}.
\begin{figure}[h]
\centerline{\includegraphics[scale=.8]{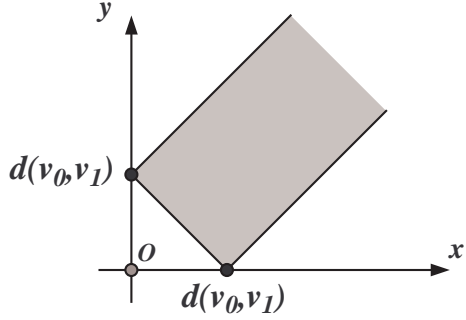} }
\caption{Area of possible triangle side lengths.}
\label{fi:triangle_area}
\end{figure}
The function of two variables $\frac{g(x,y)}{d(v_0,v_2)}$ attains
the minimal value $d^2(v_0,v_1)+ d(v_0,v_1) d(v_0,v_2)$ on the
boundary of the specified area. Hence, the inequality
    $g(s) \ge d(v_0,v_2) d(v_0,v_1) (d(v_0,v_1) + d(v_0,v_2) )$
holds for any $s$ in the component of $v_1$.

{\sc Case 3.} If $s$ belongs to the component of $v_2$ then using
same arguments as for the previous case one shows that
(\ref{eq:cut_point_ineq}) holds.
\end{proof}

\begin{corollary} \label{co:cut_point2}
Let $\Gamma$ be a connected graph, $v_0$ a cut-point in $\Gamma$,
and $v_1,v_2$ belong to distinct components of $\Gamma \setminus
\{v_0\}$. Then the inequality
    $$
    \begin{array}{c}
    d(v_0,v_2) \rb{M(v_1)-M(v_0)} + d(v_0,v_1) \rb{M(v_2)-M(v_0)}\ge C > 0 \\
    \end{array}
    $$
holds, where $C = C(v_0, v_1, v_2) = d(v_0,v_2) d(v_0,v_1)
(d(v_0,v_1) + d(v_0,v_2))$.
\end{corollary}

\begin{proof}
Indeed,
    $$d(v_0,v_2) \rb{M(v_1)-M(v_0)} + d(v_0,v_1) \rb{M(v_2)-M(v_0)}$$
    $$ = \sum_{s\in V(\Gamma)} \rb{d(v_0,v_2) \rb{d^2(v_1,s)-d^2(v_0,s)} + d(v_0,v_1) \rb{d^2(v_2,s)-d^2(v_0,s)}} \mu(s)$$
    $$ \ge \sum_{s\in V(\Gamma)} C \mu(s) = C = d(v_0,v_2) d(v_0,v_1) (d(v_0,v_1) + d(v_0,v_2) ) \mbox{, by Lemma \ref{le:cut_point_ineq1}}$$
\end{proof}

\begin{corollary}[Cut Point Lemma] \label{le:M_cutpoint}
Let $\Gamma$ be a connected graph, $v_0$ a cut-point in $\Gamma$. If
$v_1$ and $v_2$ belong to distinct connected components of $\Gamma
\setminus \{v_0\}$, then the inequalities $M(v_0) \ge M(v_1)$ and
$M(v_0) \ge M(v_2)$ cannot hold simultaneously.
\end{corollary}

\begin{proof}
Assume to the contrary that  $M(v_0) \ge M(v_1)$ and $M(v_0) \ge
M(v_2)$ hold simultaneously which is equivalent to $M(v_1) - M(v_0)
\le 0$ and $M(v_2) - M(v_0) \le 0.$ Then, multiplying by positive
constants and adding the inequalities above, we get
$d(v_0,v_2)(M(v_1) - M(v_0)) + d(v_0,v_1)(M(v_2) - M(v_0)) \le 0$
which is impossible by Corollary \ref{co:cut_point2}. This
contradiction finishes the proof.
\end{proof}

\begin{corollary}[Mean-set in a graph with a cut-point]
Let $v_0$ be a cut-point in a graph $\Gamma$ and $\Gamma \setminus
\{v_0\}$ a disjoint union of connected components
$\Gamma_1,\ldots,\Gamma_k$. Then for any distribution $\mu$ on
$\Gamma$ there exists a unique $i=1,\ldots,k$ such that $\ME \mu
\subseteq V(\Gamma_i) \cup \{v_0\}$.
\end{corollary}

\begin{corollary}[Mean-set in a graph with several cut-points]
Let $v_1,\ldots,v_n$ be cut-points in a graph $\Gamma$ and $\Gamma
\setminus \{v_1,\ldots,v_n\}$ a disjoint union of connected
components $\Gamma_1,\ldots,\Gamma_k$. Then for any distribution
$\mu$ on $\Gamma$ there exists a unique $i=1,\ldots,k$ such that
$\ME \mu \subseteq V(\Gamma_i) \cup \{v_1,\ldots,v_n\}$.
\end{corollary}

\begin{corollary}
Let $G_1$ and $G_2$ be finitely generated groups and $G = G_1 \ast
G_2$ a free product of $G_1$ and $G_2$. Then for any distribution
$\mu$ on $G$ the set $\ME \mu$ is a subset of elements of the forms
$g G_1$ or $g G_2$ for some element $g \in G$.
\end{corollary}

\begin{proposition}\label{pr:tree_centers}
Let $\Gamma$ be a tree and $\mu$ a probability measure on
$V(\Gamma)$. Then $|\ME \mu| \le 2$. Moreover, if $\ME \mu =\{u,v\}$
then $u$ and $v$ are adjacent in $\Gamma$.
\end{proposition}

\begin{proof}
Observe that any points $v_1,v_0,v_2$ such that $v_0$ is connected
to $v_1$ and $v_2$ satisfy the assumptions of Cut Point Lemma
(Corollary \ref{le:M_cutpoint}). Assume that $v_0 \in \ME \mu$. At
most one of the the neighbors of $v_0$ can belong to $\ME \mu$,
otherwise we would have $3$ connected vertices with equal $M$ values
which contradicts Cut Point Lemma.

\end{proof}

\begin{corollary}
Let $\mu$ be a probability distribution on a free group $F$. Then
$|\ME \mu| \le 2$.
\end{corollary}

In general, the number of central points can be unlimited. To see
this, consider the complete graph $K_n$ on $n$ vertices and let
$\mu$ be a uniform probability distribution on $V(K_n)$. Clearly
$\ME \mu = V(K_n)$. Another example of the same type is a cyclic
graph $C_n$ on $n$ vertices with a uniform probability distribution
$\mu$ on $V(C_n)$. Clearly $\ME \mu = V(C_n)$.
In all previous examples, the centers in a graph formed a connected
subgraph. This is not always the case. One can construct graphs with
as many centers as required and property that distances between
centers are very large (as large as one wishes).

\section{Computation of mean-sets in graphs}\label{se: Computation_meanset}

In this section we discuss computational issues that we face in
practice.
One of the technical difficulties is that, unlike the average value
$S_n/n$ for real-valued random variables, the sample mean-set $\MS_n
\equiv \MS(\xi_1,\ldots, \xi_n)$ is hard to compute. Let $G$ be a
group and $\{\xi\}_{i=1}^n$ a sequence of random i.i.d. elements
taking values in $G$. Several problems arise when trying to compute
$\MS_n$:
\begin{itemize}
    \item
Computation of the set $\{M(g) \mid g\in G\}$ requires $O(|G|^2)$
steps. This is computationally infeasible for large $G$, and simply
impossible for infinite groups. Hence  we might want to reduce the
search of a minimum to some small part of $G$.
    \item
There exist infinite groups in which the distance function $|\cdot|$
is very difficult to compute. The braid group $B_\infty$ is one of
such groups. The computation of the distance function for $B_\infty$
is an NP-hard problem, see \cite{PRaz}. Such groups require special
treatment. Moreover, there exist infinite groups for which the
distance function $|\cdot|$ is not computable. We omit consideration
of such groups.
\end{itemize}
On the other hand, we can try to devise some heuristic procedure for
this task. As we show below, if the function $M$ satisfies certain
local monotonicity properties, then we can achieve good results. The
next algorithm is a simple direct descent heuristic which can be
used to compute the minimum of a function $f$.

\begin{algorithm}{\bf (Direct Descent Heuristic)} \label{al:direct_descend}
\\{\sc Input:} A graph $\Gamma$ and a function $f:V(\Gamma) \rightarrow \MR$.
\\{\sc Output:} A vertex $v$ that locally minimizes $f$ on $\Gamma$.
\\{\sc Computations:}
\begin{itemize}
    \item[A.]
Choose a random $v \in V(\Gamma)$.
    \item[B.]
If $v$ has no adjacent vertex with smaller value of $f$, then output
current $v$.
    \item[C.]
Otherwise put $v \leftarrow u$ where $u$ is any adjacent vertex such
that $f(u) < f(v)$. Go to step B.
\end{itemize}
\end{algorithm}
It turns out that if a function $f$ satisfies certain local
properties, then we can achieve good results; namely, the proposed
algorithm finds the vertex that minimizes $f$ on $\Gamma$ exactly.
We say that a function $f:V(\Gamma) \rightarrow \MR$ is {\em locally
decreasing} if at any vertex $v\in V(\Gamma)$, such that $f$ does
not have minimum at $v$, there exists an adjacent vertex $u$ such
that $f(u) < f(v)$. We say that a function $f$ is {\em locally
finite} if for any $a,b \in \MR$ the set $f(V(\Gamma)) \cap [a,b]$
is finite.

\begin{lemma}\label{le:descend_for_convex}
Let $\Gamma$ be a graph and $f:V(\Gamma) \rightarrow \MR$ a
real-valued function that attains its minimum on $\Gamma$. If $f$ is
locally decreasing and locally finite, then Algorithm
\ref{al:direct_descend} for $\Gamma$ and $f$ finds the vertex that
minimizes $f$ on $\Gamma$.
\end{lemma}

\begin{proof}
Let $v \in V(\Gamma)$ be a random vertex chosen by Algorithm
\ref{al:direct_descend} at Step A. If $v$ is a minimum of $f$, then
the algorithm stops with the correct answer $v$. Otherwise, the
algorithm, at Step C, chooses any vertex $u$ adjacent to $v$ such
that $f(u) < f(v)$. Such a vertex $u$ exists, since the function $f$
is locally decreasing by assumption. Next, Algorithm
\ref{al:direct_descend} performs the same steps for $u$.
Essentially, it produces a succession of vertices $v_0,v_1,v_2,
\ldots $ such that $v_0 = v$ and, for every $i= 0, 1, 2, \ldots$,
the vertices $v_i,v_{i+1}$ are adjacent in $\Gamma$ with the
property $f(v_i)>f(v_{i+1})$.

We claim that the constructed succession cannot be infinite. Assume,
to the contrary, that the chain $v_0,v_1,v_2, \ldots $ is infinite.
Let $m$ be the minimal value of $f$ on $\Gamma$. Then $f(V(\Gamma))
\cap [m,f(v)]$ is infinite, and, $f$ cannot be locally finite.
Contradiction. Hence the sequence is finite, and the last vertex
minimizes $f$ on $V(\Gamma)$.
\end{proof}

\begin{lemma}\label{le:tree_loc_finite}
Let $\mu$ be a distribution on a locally finite graph $\Gamma$ such
that a weight function $M(\cdot)$ is defined. Then the function
$M(\cdot)$ is locally finite on $\Gamma$.
\end{lemma}

\begin{proof}
Since the function $M$ is non-negative, it suffices to prove that
for any $b\in \MR_+$ the set $M(V(\Gamma)) \cap [0,b]$ is finite.
Let $v \in \ME(\xi)$, i.e., $v$ minimizes the value of $M$, and
$r\in\MN$ such that
    $0 < \frac{1}{2} M(v) \le \sum_{i\in B_v(r)} d^2(v,i) \mu(i),$
as in the proof of Lemma \ref{le:E_finite}. Choose an arbitrary
value $b\in \MR_+$ and put $\alpha = \max\{2,b/M(v)\}$. Then one can
prove (as in Lemma \ref{le:E_finite}) that for any $u\in \Gamma
\setminus B_v((\alpha+2)r)$, we have
    $M(u)> (\alpha+1) M(v) > b.$
Therefore, $M(V(\Gamma)) \cap [0,b] \subset M(B_v((\alpha+2)r))$ and
the set $B_v((\alpha+2)r)$ is finite.
\end{proof}

\begin{theorem} \label{thm:Algorithm_works_trees}
Let $\mu$ be a distribution on a locally finite tree $T$ such that a
function $M$ is totally defined. Then Algorithm
\ref{al:direct_descend} for $T$ and $M$ finds a central point
(mean-set) of $\mu$ on $T$.
\end{theorem}

\begin{proof}
Follows from Lemmata \ref{le:descend_for_convex},
\ref{le:M_cutpoint}, \ref{le:tree_loc_finite}, and
\ref{le:E_finite}.
\end{proof}

Note, the function $M$ is not locally decreasing for every
graph, and a local minimum, computed by Algorithm \ref{al:direct_descend},
is not always a global minimum.

\section{Experiments} \label{se: experiments}

In this section we demonstrate how the technique of
computing mean-sets, employing the Direct Descent Algorithm
\ref{al:direct_descend} described in  section \ref{se:
Computation_meanset}, works in practice and produces results
supporting our SLLN for graphs and groups. More precisely, we
arrange series of experiments in which we compute the sample
mean-sets of randomly generated samples of $n$ random elements and
observe a universal phenomenon: the greater the sample size $n$, the
closer the sample mean gets to the actual mean of a given
distribution. In particular, we experiment with free groups,
in which the length function is easily
computable. All experiments are done using the CRAG software
package, see \cite{CRAG}.

One of the most frequently used distributions on the free groups is
a uniform distribution $\mu_L$ on a {\em sphere} of radius $L$
defined as
$S_L = \{w\in F(X) \mid |w|=L\}.$
Clearly, $S_L$ is finite. Therefore, we can easily define a uniform
distribution $\mu_L$ on it as
$$
\mu_L(w) = \left\{
\begin{array}{ll}
\frac{1}{|S_L|} & \mbox{if } |w|=L;\\
0 & \mbox{otherwise}.\\
\end{array}
\right.
$$

The reader interested in the question of defining probabilities on groups can find
several approaches to this issue in
\cite{BoMS}.
One of the properties of $\mu_L$ is that its mean-set is just the
trivial element of the free group $F(X)$.
Observe also that the distance of any element of
$F(X)$ to the mean-set is just the length of this element (or
length of the corresponding word).

Table \ref{tb:free_group4} below contains
the results of experiments for the distributions
$\mu_5$, $\mu_{10}$, $\mu_{20}$, $\mu_{50}$ on the group $F_4$.
The main parameters in our experiments are the rank $r$ of the free
group, the length $L$, and the sample size $n$. For every particular
triple of parameter values $(r,L,n)$, we perform series of $1000$
experiments to which we refer (in what follows), somewhat loosely,
as series $(r,L,n)$. Each cell in the tables below corresponds to a
certain series of experiments with parameters $(r,L,n)$. In each
experiment from the series $(r,L,n)$, we randomly generate $n$ words
$w_1,\ldots,w_n$, according to distribution $\mu_L$, compute the
sample mean-set $\MS_n$, and compute the displacement
of the actual center of $\mu_L$ from
$\MS_n$. The set $\MS_n$ is computed
using Algorithm \ref{al:direct_descend} which, according to Theorem
\ref{thm:Algorithm_works_trees}, always produces correct answers for
free groups. Every cell in the tables below contains a pair of
numbers $(d,N)$; it means that in $N$ experiments out of $1000$ the
displacement from the real mean was $d$.

  \begin{table}[ht]
  \tiny
  \centering
  \begin{tabular*}{1\textwidth}%
     {@{\extracolsep{\fill}}|c||c|c|c|c|c|c|c|c|}
 \hline
 {\bf L$\backslash$n} &{\bf 2} & {\bf 4} & {\bf 6} & {\bf 8} & {\bf 10} & {\bf 12} & {\bf 14} & {\bf 16} \\
 \hline
 \hline
{$\bf\mu_{5}$}  & (0,885) & (0,943) & (0,978) & (0,988) & (0,999) & (0,998) & (0,1000) & (0,999) \\
 & (1,101) & (1,55) & (1,22) & (1,12) & (1,1) & (1,2) &  & (1,1) \\
 & (2,13) & (2,2) &  &  &  &  &  &  \\
 & (3,1) &  &  &  &  &  &  &  \\
\hline
{$\bf\mu_{10}$} & (0,864) & (0,930) & (0,976) & (0,993) & (0,994) & (0,999) & (0,1000) & (0,1000) \\
 & (1,117) & (1,69) & (1,24) & (1,7) & (1,6) & (1,1) &  & \\
 & (2,16) & (2,1) &  &  &  &  &  &  \\
 & (3,2) &  &  &  &  &  &  &  \\
 & (4,1) &  &  &  &  &  &  &  \\
\hline
{$\bf\mu_{20}$} & (0,859) & (0,940) & (0,975) & (0,985) & (0,991) & (0,1000) & (0,999) & (0,999) \\
 & (1,116) & (1,58) & (1,25) & (1,15) & (1,9) &  & (1,1) & (1,1) \\
 & (2,19) & (2,2) &  &  &  &  &  &  \\
 & (3,6) &  &  &  &  &  &  &  \\
\hline
{$\bf\mu_{50}$} & (0,872) & (0,928) & (0,984) & (0,991) & (0,998) & (0,997) & (0,998) & (0,999) \\
 & (1,108) & (1,71) & (1,16) & (1,9) & (1,2) & (1,3) & (1,2) & (1,1)  \\
 & (2,19) & (2,1) &  &  &  &  &  &  \\
 & (3,1) &  &  &  &  &  &  &   \\
\hline
  \end{tabular*}
  \caption{The results of experiment for $F_4$.}
  \label{tb:free_group4}
  \end{table}

By doing experiments for free groups of higher ranks, one can easily observe that
as the rank of the free
group grows, we get better and faster convergence.
Intuitively, one
may think about this outcome as follows: the greater the rank is,
the more branching in the corresponding Cayley graph we have, which
means that more elements are concentrated in a ball, and the bigger
growth (in that sense) causes the better and faster convergence.

\vspace{5mm} {\bf Acknowledgements.} We are grateful to Ioannis
Karatzas for his support and his time spent reading this manuscript.
The authors also extend their gratitude to Persi Diaconis for his suggestions of
useful literature sources. In addition, we are thankful to Gerard
Ben Arous for a discussion about n-dimensional random walks.


\providecommand{\bysame}{\leavevmode\hbox to3em{\hrulefill}\thinspace}
\providecommand{\MR}{\relax\ifhmode\unskip\space\fi MR }
\providecommand{\MRhref}[2]{%
  \href{http://www.ams.org/mathscinet-getitem?mr=#1}{#2}
}
\providecommand{\href}[2]{#2}

\end{document}